\documentclass[twoside,11pt,leqno]{article}

\usepackage{amsmath,amssymb}
\usepackage{enumerate}

\let\bsbf=\bf

\newtheorem{theorem}{Theorem}[section]
\newtheorem{lemma}[theorem]{Lemma}
\newtheorem{corollary}[theorem]{Corollary}

\newtheorem{remark}[theorem]{Remark}
\newtheorem{example}[theorem]{Example}
\newtheorem{definition}[theorem]{Definition}
\newtheorem{conjecture}[theorem]{Conjecture}

\newcommand{\rnf}{\renewcommand{\thefootnote}{\arabic{footnote}}}

\newcommand{\aire}{\vphantom{A^{A^A}}}
\newcommand{\maire}{\vphantom{A^{A^{A^a}}}}
\newcommand{\mmaire}{\vphantom{A^{A^{A^A}}}}

\title{\vspace{-3cm}
\ \hspace{-2.70in}
{\small {\bf Morfismos},
Vol. 13, No. 2, 2009,
pp.
        1--\pageref{ultimapagina} 
}
\\
\vspace{2.6cm} 
The Arf-Kervaire Invariant of framed manifolds
\footnote{\bsbf Invited article.}
}

\author{\rnf   Victor P.~Snaith
}

\date{}

\begin{document}
\maketitle






\begin{abstract} \noindent
This work surveys classical and recent advances around the existence
of exotic differentiable structures on spheres and its connection 
to stable homotopy theory.
\end{abstract}

\noindent 


\vspace{.7cm}\centerline{\sc Contents}
1 ~\hspace{.6mm}Introduction\hfill\pageref{introduccion}\ \ \ \ 

2 ~\hspace{.6mm}Stable homotopy groups of spheres\hfill\pageref{esferas}\ \ \ \ 

3 ~\hspace{.6mm}Framed manifolds and stable homotopy groups\hfill\pageref{enmarcados}\ \ \ \ 

4 ~\hspace{.6mm}The classical stable homotopy category\hfill\pageref{categoria}\ \ \ \ 

5 ~\hspace{.6mm}Cohomology operations\hfill\pageref{operaciones}\ \ \ \ 

6 ~\hspace{.6mm}The classical Adams spectral sequence\hfill\pageref{ASS}\ \ \ \ 

7 ~\hspace{.6mm}Homology operations\hfill\pageref{homologia}\ \ \ \ 

8 ~\hspace{.6mm}The Arf-Kervaire invariant one problem\hfill\pageref{AKi1}\ \ \ \ 

9 ~\hspace{.6mm}Exotic spheres, the J-homomorphism and the\hfill \ \ \ \ \ 

~~~~Arf-Kervaire invariant\hfill\pageref{Jmorfismo}\ \ \ \ 

\ \hspace{-2.5mm}10 ~\hspace{.1mm}Non-existence results for the Arf-Kervaire invariant\hfill
\pageref{noexiste}\ \ \ \ 

\ \hspace{-2.5mm}References\hfill\pageref{referencias}\ \ \ \ 

\vspace{.22cm}
\begin{quotation}
``I know what you're thinking about,'' said Tweedledum; ``but it isn't so, nohow.'' ``Contrariwise,''  continued Tweedledee,``if it was so, it might be; and if it were so, it would be; but as it isn't; it ain't.  That's logic.''
\vspace{5pt}

{\em Through the Looking Glass}
by Lewis Carroll  (aka Charles Lutwidge Dodgson)  \cite{Car40}
\vspace{5pt}

\end{quotation}

\section{Introduction}\label{introduccion}
 This survey article deals with significant recent progress concerning the existence of continuous maps $\theta : S^{N+m} \longrightarrow S^{N}$ between high-dimensional spheres having non-zero Arf-Kervaire invariant. Even to sketch the problem will take some preparation, at the moment it will suffice to say that this is a fundamental unresolved problem in homotopy theory which has a history extending back over fifty years. During that time such maps have been constructed for only five values of $m$. To a homotopy theorist it qualifies as a contender for the most important unsolved problem in $2$-adic stable homotopy theory. 
 
\vspace{.7mm}
I was invited to write this survey article for two reasons. Firstly because I had recently published a monograph on the Arf-Kervaire invariant \cite{Arfbk08} in which I conjectured what the experts must have been feeling for years, that such $\theta$'s only existed for the values $m = 2, 6, 14, 30, 62$. Secondly, much more importantly, because recently preprints have appeared (\cite{Akh08a}  and \cite{HHR09}) which almost completely prove this non-existence 
conjecture---only $m = 126$ remains unresolved.
 
\smallskip
Since the inception of algebraic topology \cite{Poin1895} the study of homotopy classes of continuous maps between spheres has enjoyed a very exceptional, central role. As is well-known, for homotopy classes of maps $f: S^{n} \longrightarrow  S^{n}$ with $n \geq 1$ the sole homotopy invariant is the degree, which characterises the homotopy class completely. The search for a continuous map between spheres of different dimensions and not homotopic to the constant map had to wait for its resolution until the remarkable paper of Heinz\index{Heinz Hopf} Hopf  \cite{Hopf35}. In retrospect, finding an example was rather easy because there is a canonical quotient map from $S^{3} $ to the orbit space of the free circle action 

\vspace{-2mm}
$$S^{3}/S^{1} = {\mathbb CP}^{1}= S^{2}.$$ 

\vspace{-1mm}
\noindent On the other hand, the problem of showing that this map is not homotopic to the constant map requires either ingenuity (in this case Hopf's observation that the inverse images of any two distinct points on $S^{2}$ are linked circles) or, more influentially, an invariant which does the job (in this case the Hopf invariant). The Hopf invariant is an integer which is associated to any continuous map of the form $f: S^{2n-1} \longrightarrow  S^{n}$ for $n \geq 1$. Hopf showed that when $n$ is even there exists a continuous map whose Hopf invariant is equal to any even integer. On the other hand the homotopy classes of continuous maps  $g:  S^{m} \longrightarrow  S^{n}$ in almost all cases with 
$m > n \geq 1$ form a finite abelian group. For the study of the $2$-Sylow subgroup\index{Sylow subgroup} of these groups the appropriate invariant is the Hopf invariant modulo $2$. With the construction of mod $p$ cohomology operations by Norman Steenrod\index{Norman Steenrod} it became possible to define the mod $2$ Hopf invariant for any $g$ but the only possibilities for non-zero mod $2$ Hopf invariants occur when $m-n+1$ is a power of two (\cite[p.~12]{StE62}). 

\vspace{.1mm}
When $n \gg 0$ the homotopy classes of such $g$'s form the stable homotopy group $\pi_{m-n}( \Sigma^{\infty} S^{0})$, (spectra and their homotopy groups will be explained in \S\ref{1.3.1} and Example \ref{1.3.2})
which is a finite abelian group when $m > n$.
The $p$-Sylow subgroups of stable homotopy groups were first organised systematically by the 
mod $p$ Adams spectral sequence, constructed by Frank Adams in \cite{Ad58}.
Historically, the case when $p=2$ predominates. Spectral sequences require a few words of explanation, after which we shall see in Theorem \ref{1.1.2} that, on the line $s=1$ in the classical Adams spectral sequence four elements exist denoted by $h_{0}, h_{1}, h_{2}, h_{3}$. They correspond, respectively, to four stable homotopy elements in $\pi_{j}( \Sigma^{\infty} S^{0})$ when $j=0, 1, 3$ and $7$, respectively. In positive dimensions the 
homotopy classes with non-zero mod $2$ Hopf invariant would all be represented on the $s=1$ line in dimensions of the form $j = 2^{k}-1$, by Steenrod's result. However, a famous result due originally to Frank Adams (\cite{Ad60}; see also \cite{AdAt66}, \cite{Sn83} and \cite[Theorem~6.3.2]{Arfbk08}) shows that only $h_{1}, h_{2}, h_{3}$ actually correspond to homotopy classes with non-zero mod $2$ Hopf invariant.

The non-existence of  homotopy classes with Hopf invariant one was a big-ticket item in its day. Inspection of the segment of the spectral sequence which is given in Theorem \ref{1.1.2} correctly gives the impression that next important problem concerns whether or not the classes labelled $h_{i}^{2}$ represent elements of $\pi_{2^{i+1}-2}( \Sigma^{\infty} S^{0})$. The invariant which is capable of detecting homotopy classes represented on the $s=2$ line is due to Michel Kervaire \cite{Kerv60} as generalised by Ed Brown Jr.~\cite{EHB72}. Bill Browder discovered the fundamental result \cite{Brow69}
(this is Theorem \ref{9.3}, whose proof by the technique of \cite{Arfbk08} is sketched in \ref{9.4}), the analogue of Steenrod's result about the Hopf invariant, that the Arf-Kervaire invariant could only detect stable homotopy classes in dimensions of the form $2^{i+1}-2$.

\vspace{.11mm}I am very grateful to Peter Landweber for a series of email discussions during 2008 concerning \cite{Akh08a} and \cite{Akh08}. Similarly I am very grateful to Mike Hopkins and Doug Ravenel for an email and a phone call  respectively in 2009\hspace{.2mm}---a few days before Mike announced the result of \cite{HHR09} at the Atiyah Fest in Edinburgh, Scotland and Doug lectured on it at a conference in Lisbon, Portugal. In addition, I am immensely indebted to the referees for their meticulous reading of this review.

\section{Stable homotopy groups of spheres}\label{esferas}
\subsection{Digression on spectral sequences}
\label{1.1.0}

In a textbook dealing with spectral sequences (e.g.~\cite[Chapter~VIII]{HS71},
\cite{Mc01}, or \cite[Chapter~9]{Spa66}) 
one is apt to find a notation of the form
$   E_{p,q}^{2}   \Rightarrow  E_{p+q} .   $
These symbols denote an organisational structure called a spectral sequence and useful for computing the graded group $E_{*} = \bigoplus_{n \in {\mathbb Z}} \ E_{n}$. Notice that $E^{2}_{*,*} =   \bigoplus_{p,q} \   E_{p,q}^{2}$ is a bi-graded group
and bi-gradings in spectral sequences can take a variety of forms. However, typically, there are families of {\em differentials} which are homomorphisms of the form $d_{r} :  E_{p,q}^{r}  \longrightarrow  E_{p-r,q+r-1}^{r}$ for each $r \geq 2$ and satisfying (i) $d_{r} \cdot d_{r} = 0$ and (ii) there are given isomorphisms of the form
\[   E_{p,q}^{r+1}  \cong  \frac{{\rm Ker}(d_{r} :  E_{p,q}^{r}  \longrightarrow  E_{p-r,q+r-1}^{r})}{{\rm Im}(d_{r} :  E_{p+r,q-r+1}^{r}  \longrightarrow  E_{p,q}^{r})}  .  \]
Knowledge of all the differentials enables one to compute the sequence of groups
$  E_{p,q}^{2},  E_{p,q}^{3}, \ldots ,   E_{p,q}^{r}, \ldots $. In useful cases this sequence asymptotes to a common value, which is a group denoted by $E_{p,q}^{\infty}$ of which a particularly handy example occurs when  
$E_{p,q}^{\infty} = 0$.
Given all this, the final part of the structure of a spectral sequence consists of a filtration on each $E_{m}$ of the form
$  \cdots   \subseteq  F_{2}E_{m}  \subseteq     F_{1}E_{m}  \subseteq  F_{0}E_{m} = E_{m} $ together with isomorphisms of the form
$  F_{p}E_{p+q}/F_{p+1}E_{p+q} \cong   E_{p,q}^{\infty}$.

\smallskip
In the event that an element $x \in\aire E_{p,q}^{2} $ satisfies $d_{2}(x)=0$ and that its homology class in $E_{p,q}^{3} $
lies in the kernel of $d_{3}$ and so on all the way to $E_{p,q}^{\infty} $ we say that $x$ is an {\em infinite cycle}. Furthermore if $z \in E_{p+q}$ lies in $F_{p}E_{p+q}$ and coincides in $  F_{p}E_{p+q}/F_{p+1}E_{p+q} $ with the coset associated with the infinite cycle $x$ we say that $z$ is {\em represented by} $x$ in the spectral sequence. Usually, in the presence of lots of non-zero differentials, the relation between $z$ and $x$ has a lot of indeterminacy. However, if all the $d_{r}$ vanish, which is referred to as the {\em collapsing} at  $E^{2}_{*,*} $ of the spectral sequence, then one can deduce a lot about $z$ (e.g. its order) from properties of $x$.

\smallskip
In the case of the classical $2$-adic Adams spectral sequence, $E_{m} = \pi_{r}^{S}(S^{0}) \otimes {\mathbb Z}_{2}$ where ${\mathbb Z}_{2}$ denotes the $2$-adic integers and $\pi_{r}^{S}(S^{0})$ is the stable stem defined in \S\ref{1.1.1}.

\subsection{The History}
\label{1.1.1}
The set of homotopy classes of continuous maps from $S^{m}$ to $S^{n}$ form a group denoted by $\pi_{m}(S^{n})$. To add two maps we squeeze the equator of $S^{m}$ to a point, making two copies of $S^{m}$ attached at one point, and we put one of the given maps on each $S^{m}$. There is a homomorphism from $\pi_{m}(S^{n})$ to $\pi_{m+1}(S^{n+1})$ given by treating the original spheres as equators in $S^{m+1}$ and $S^{n+1}$ and mapping each latitudinal $S^{m}$ to the $S^{n}$ at the same latitude by a scaled down version of the original map. The stable homotopy group (sometimes called a {\em stable stem}) is the direct limit $\pi_{k}^{S}(S^{0}) = \lim_{r \rightarrow \infty} \ \pi_{r+k}(S^{r})\,$---also denoted by $\pi_{k}(\Sigma^{\infty} S^{0})$.

\smallskip
 The calculation of the stable homotopy groups of spheres is one of the most central and intractable problems in algebraic topology. In the 1950's Jean-Pierre Serre used his spectral sequence to study the problem \cite{Ser58}. In 1962 Hirosi Toda used his triple products and the EHP exact sequence to calculate the first nineteen stems (that is, $\pi_{j}^{S}(S^{0})$ for $0 \leq j \leq 19$) \cite{Toda62}. These methods were later extended by Mimura, Mori, Oda and Toda to compute the first thirty stems 
(\cite{ Mim64},  \cite{MTod63},  \cite{MMOda75},   \cite{Oda79}).  In the late 1950's the study of the classical Adams spectral sequence began \cite{Ad58}.  According to \cite[Chapter~1]{Ko90}, computations in this spectral sequence were still being pursued into the 1990's using the May spectral sequence and the lambda algebra. The best published results are Peter May's thesis (\cite{Ma64}, \cite{Ma66}) and the computation of the first forty-five stems by Michael Barratt, Mark Mahowald and Martin Tangora (\cite{BMT70}, \cite{MT67}), as corrected by Bob Bruner \cite{Brun84}. The use of the Adams\index{Frank Adams} spectral sequence based on Brown-Peterson\index{Frank Peterson} cohomology\index{Ed Brown Jr.} theory ($BP$ theory for short) was initiated by Sergei Novikov \cite{Nov67} and Raph Zahler \cite{Zah72}. The $BP$ spectral sequences were most successful at odd primes \cite{MRW77}. A comprehensive survey of these computations and the methods which have been used is to be found in Doug Ravenel's\index{Doug Ravenel} book \cite{Rav86}.

\smallskip
We shall be interested in two-primary phenomena in the stable homotopy of spheres, generally ackowledged to be the most intractable case. What can be said at the prime $p=2$? A
seemingly eccentric, cart-before-the-horse method for computing stable stems was developed in 1970 by Joel Cohen \cite{Coh70}. For a generalised homology theory $E_{*}$ and a stable homotopy spectrum $X$ there is an Atiyah-Hirzebruch spectral\index{Friedrich Hirzebruch} sequence \cite{Dyer69}
\[   E_{p,q}^{2} = H_{p}( X ; \pi_{q}(E))  \Longrightarrow  E_{p+q}(X) .   \]
Cohen studied this spectral sequence with $X$ an Eilenberg-Mac\hspace{1pt}Lane spectrum and $E$ equal to stable homotopy or stable homotopy modulo $n$. The idea was that in this case the spectral sequence is converging to zero in positive degrees and, since the homology of the Eilenberg-Mac\hspace{1pt}Lane spectra
are known, one can inductively deduce the stable homotopy groups of spheres. This strategy is analogous to the inductive computation of the cohomology of Eilenberg-Mac\hspace{1pt}Lane spaces by means of the Serre spectral sequence \cite{Cartan54}. Cohen was only able to compute a few low-dimensional stems before the method became too complicated. Therefore the method was discarded in favour of new methods which used the Adams spectral sequence. In 1972 Nigel Ray used the Cohen method with $X = MSU$ and $E = MSp$, taking advantage of knowledge of $H_{*}(MSU)$ and $MSp_{*}(MSU)$ to compute $\pi_{j}(MSp)$ for $0 \leq j \leq 19$ \cite{Ray72}. Again this method was discarded because David Segal had computed up to dimension thirty-one via the Adams spectral sequence and in \cite{Ko93} these computations were 
pushed to dimension one hundred.

\smallskip
In 1978\hspace{.5mm} Stan Kochman\index{Stan Kochman} and I studied the Atiyah-Hirzebruch spectral sequence method in the case where $X = BSp$ and $E_{*}$ is stable homotopy \cite{KSn79}.  An extra ingredient was added in this paper; namely we exploited the Landweber-Novikov\index{Peter Landweber} operations\index{Sergei Novikov} to study differentials. This improvement is shared by the case when $X = BP$ and $E_{*}$ is stable homotopy, capitalising on the sparseness of $H_{*}(BP)$ and 
using Quillen\index{Dan Quillen} operations to compute the differentials. The result of this method is a computer-assisted calculation of the stable stems up to dimension sixty-four with particular emphasis on the two-primary part  \cite{Ko90}. 

\smallskip
Let ${\mathcal A}$ denote the mod $2$ Steenrod algebra (see Definition \ref{1.6.3}). From \cite[Theorem.~3.2.11]{Rav86} we have the following result. 
\begin{theorem}
\label{1.1.2}
For the range $t-s \leq 13$ and $s \leq 7$ the group ${\rm Ext}_{{\mathcal A}}^{s,t}( {\mathbb Z}/2 , 
$ ${\mathbb Z}/2)$ is generated as an ${\mathbb F}_{2}$-vector-space by the elements listed in the following table. Only the non-zero groups and elements are tabulated, the vertical coordinate is $s$ and the horizontal is $t-s$.
\[ \begin{array}{c|cccccccccc}
7 & h_{0}^{7} & &&&&&&&& P(h_{1}^{2}) \\[1ex]
6 & h_{0}^{6} &&&&&&&& P(h_{1}^{2}) &  P(h_{0}h_{2}) \\[1ex]
5 & h_{0}^{5} &&&&&&&  P(h_{1}) && P(h_{2})  \\[1ex]
4& h_{0}^{4} &&&&& h_{0}^{3}h_{3}&& h_{1}c_{0}& &\\[1ex]
3 &h_{0}^{3}  &&& h_{1}^{3} && h_{0}^{2}h_{3}& c_{0} & h_{1}^{2}h_{3} && \\[1ex]
2& h_{0}^{2} && h_{1}^{2}& h_{0}h_{2}& h_{2}^{2}& h_{0}h_{3}& h_{1}h_{3} &&&\\[.4ex]
1& h_{0}  & h_{1}&& h_{2} && h_{3} &&&&\\
0 &  1  &&&&&&&&& \\
\hline
&0 &1&2&3&6&7&8&9&10 & 11  \\
\end{array} \]
There are no generators for $t-s = 12, 13$ and the only generators in this range with $s >7$ are powers of $h_{0}$.
\end{theorem}

Inspecting this table one sees that there can be no differentials in this range and we obtain the following table of values for the $2$-Sylow subgroups of the stable stems $\pi_{n}^{S}(S^{0}) \otimes {\mathbb Z}_{2}$.

\begin{corollary}
\label{1.1.3}

For $n \leq 13$ the non-zero groups $\pi_{n}^{S}(S^{0}) \otimes {\mathbb Z}_{2}$ are given by the following table.
\begin{center}
\begin{tabular}{|c|c|c|c|c|c|} \hline
$n$ & $0$ & $1$ & $2$ & $3$ & $6$ \\
\hline
$\pi_{n}(S^{0}) \otimes {\mathbb Z}_{2} $ & $ {\mathbb Z}_{2}$ & ${\mathbb Z}/2$ & $ {\mathbb Z}/2$ & $ {\mathbb Z}/8$ &  $ {\mathbb Z}/2$ \\
 \hline
\end{tabular}\medskip

\begin{tabular}{|c|c|c|c|c|c|} \hline
$n$ & $7$ & $8$ & $9$ & $10$ & $11$ \\
\hline
$\pi_{n}(S^{0}) \otimes {\mathbb Z}_{2} $ & $ {\mathbb Z}/16$ &  $( {\mathbb Z}/2 )^{2}$ &  $( {\mathbb Z}/2 )^{3}$ & $  {\mathbb Z}/2$ &  $ {\mathbb Z}/8 $\\
 \hline
\end{tabular}
\end{center}
\end{corollary}

For more recent computational details the reader is referred to \cite{BK94}, \cite{BK96}, \cite{Ko73}, \cite{KSn79}, \cite{Ko90}, \cite{Ko91}, \cite{Ko92}, \cite{Ko92b}, \cite{Ko93}, \cite{Ko96} and \cite{KM95}.
For example, according to \cite{Ko90}, the $2$-Sylow subgroup of $\aire\pi_{62}^{S}(S^{0})$ is ${\mathbb Z}/2 \oplus   {\mathbb Z}/2  \oplus  {\mathbb Z}/4$ with an element of Arf-Kervaire invariant one denoted by $A[62,1]$ having order two.

 \section{Framed manifolds and stable homotopy \\groups}\label{enmarcados}
\begin{definition}{\em
\label{1.2.1}
Let $M^{n}$ be a compact ${\rm C}^{\infty}$ manifold\index{compact manifold} without boundary and let $i : M^{n}  \longrightarrow  {\mathbb R}^{n+r}$ be an embedding. The normal bundle\index{normal bundle} of $i$, denoted by $\nu(M, i)$, is the quotient of the pullback\index{pullback} of the tangent bundle\index{ tangent bundle}
of ${\mathbb R}^{n+r}$ by the sub-bundle given by the tangent bundle of $M$
$$ \nu(M, i)   =  \frac{  i^{*} \tau({\mathbb R}^{n+r})}{  \tau(M)}  $$

\vspace{-2mm}\noindent
so that $\nu(M, i)$ is an $r$-dimensional real vector bundle over $M^{n}$. If we give 
$\tau({\mathbb R}^{n+r}) = {\mathbb R}^{n+r}   \times  {\mathbb R}^{n+r}$ the Riemannian metric\index{Riemannian metric} obtained from the usual inner product in Euclidean space, the
total space of the normal bundle $\nu(M, i)$ may be identified with the orthogonal complement of
$\tau(M)$ in $i^{*} \tau({\mathbb R}^{n+r})$. That is, the fibre at $z \in M$ may be identified with the subspace of vectors $(i(z),x) \in  {\mathbb R}^{n+r} \times {\mathbb R}^{n+r}$ such that $x$ is orthogonal to
$i_{*} \tau(M)_{z}$ where $i_{*}$ is the induced embedding of $\tau(M)$ into $\tau({\mathbb R}^{n+r})$.
}\end{definition}
\begin{lemma}
\label{1.2.2}

If $r$ is sufficiently large (depending only on $n$) and $i_{1}, i_{2} : M^{n} \longrightarrow  {\mathbb R}^{n+r}$ are two embeddings then $\nu(M, i_{1} )$ is trivial (i.e. $\nu(M, i_{1})$ $  \cong M \times 
 {\mathbb R}^{r}$) if and only if $\nu(M, i_{2} )$ is trivial. More generally, for large $r$, the two normal bundles are isomorphic.
\end{lemma}

\begin{definition}{\em
\label{1.2.3}
Let $\xi$ be a vector bundle over a compact manifold $M$ endowed with a Riemannian metric
on the fibres. Then the Thom\index{Ren\'{e} Thom} space\index{Thom space} is defined to be the quotient of the unit disc bundle $D(\xi)$ of $\xi$ with the unit sphere bundle $S(\xi)$ collapsed to a point. Hence

\vspace{-1.8mm}\ 
$$ T(\xi)  =  \frac{D(\xi) }{ S(\xi)}  $$
is a compact topological space with a basepoint\index{basepoint} given by the image of $S(\xi)$.

\smallskip
If $M$ admits an embedding with a trivial normal bundle, as in Lemma \ref{1.2.2}, we say that $M$ has a stably trivial normal bundle\index{stably trivial normal bundle}.  Write $M_{+}$ for the disjoint union of $M$ and a disjoint base
point. Then there is a canonical homeomorphism

\vspace{-1.8mm}\ 
$$  T(  M \times  {\mathbb R}^{r}  )  \cong   \Sigma^{r}(M_{+})  $$
between the Thom space of the trivial $r$-dimensional vector bundle and the $r$-fold suspension\index{$r$-fold suspension}
of $M_{+}$,  $(S^{r} \times  (M_{+}) )/  (S^{r} \vee  (M_{+}))  =  S^{r} \wedge  (M_{+})$.
}\end{definition}
\subsection{The Pontrjagin-Thom construction}
\label{1.2.4}

Suppose that $M^{n}$ is a manifold as in Definition \ref{1.2.1} together with a choice of trivialisation 
of normal bundle $\nu(M, i)$. This choice gives a choice of homeomorphism

\vspace{-1.8mm}\ 
$$   T( \nu(M, i))   \cong   \Sigma^{r}(M_{+}).    $$
Such a homeomorhism is called a framing\index{ framing} of $(M, i)$.
Now consider the embedding $i : M^{n}  \longrightarrow  {\mathbb R}^{n+r}$ and identify the
$(n+r)$-dimensional sphere $S^{n+r}$ with the one-point compactification $ {\mathbb R}^{n+r}  \cup  \{  \infty \}$. The Pontrjagin-Thom construction\index{Pontrjagin-Thom construction} is the map

\vspace{-1.8mm}\ 
$$    S^{n+r}  \longrightarrow    T( \nu(M, i))  $$
given by collapsing the complement of the interior of the unit disc bundle $D( \nu(M, i) )$ to the point corresponding to $S( \nu(M, i) )$ and by mapping each point of $D( \nu(M, i) )$ to itself. 

\smallskip
Identifying the $r$-dimensional sphere with the $r$-fold suspension $\Sigma^{r} S^{0}$ of the zero-dimensional sphere (i.e. two points, one the basepoint) the map which collapses $M$ to the non-basepoint yields a basepoint preserving map $ \Sigma^{r}(M_{+}) \longrightarrow  S^{r}$.

Therefore, starting from a framed manifold $M^{n}$, the Pontrjagin-Thom construction yields a based map
$$  S^{n+r}  \longrightarrow    T( \nu(M, i)) \cong    \Sigma^{r}(M_{+}) \longrightarrow  S^{r}  ,$$

\vspace{-1mm}\noindent 
whose homotopy class defines an element of $\pi_{n+r}(S^{r})$.

\subsection{Framed Cobordism}
\label{1.2.5}
Two compact, smooth $n$-dimensional manifolds without boundary are called cobordant if their disjoint union is the boundary of some compact $(n+1)$-dimensional manifold. The first description of this equivalence relation appears in Henri Poincar\'{e}'s work \cite[\S5]{Poin1895}. This paper of Poincar\'{e} is universally considered to be the origin of algebraic topology. His concept of homology is basically the same as that of cobordism as later developed by Lev Pontrjagin~\cite{Pont47,Pont55}, Ren\'{e} Thom~\cite{Thom54}, and others~\cite{Sto68}.

\smallskip
The notion of cobordism can be applied to manifolds with a specific additional structure on their stable normal bundle ($\nu(M^{n}, i)$ for an embedding of codimension $r \gg n$). In particular one can define the equivalence relation of framed cobordism\index{framed cobordism} between $n$-dimensional manifolds with chosen framings of their stable normal bundle. This yields a graded ring of framed cobordism classes $\Omega_{*}^{{\rm fr}}$ where the elements in  $\Omega_{n}^{{\rm fr}}$ are equivalence classes of compact framed $n$-manifolds without boundary. The sum is induced by disjoint union and the ring multiplication by cartesian product of manifolds.

\smallskip
More generally one may extend the framed cobordism relation to maps of the form
$f : M^{n} \longrightarrow  X$ where $M^{n}$ is a compact framed manifold without boundary and
$X$ is a fixed topological space. From such a map the Pontrjagin-Thom construction yields 

\vspace{-2mm}\ 
$$   S^{n+r}  \longrightarrow    T( \nu(M, i)) \cong    \Sigma^{r}(M_{+}) \longrightarrow   \Sigma^{r}(X_{+})  $$
whose homotopy class defines an element of $\pi_{n+r}( \Sigma^{r}(X_{+}))$.

\smallskip
The cobordism of maps $f$ yields a graded ring of framed cobordism classes $\Omega_{*}^{{\rm fr}}(X)$.
The Pontrjagin-Thom element associated to $f$ does not depend solely on its class in 
$\Omega_{n}^{{\rm fr}}(X)$ but the image of $f$ in the stable homotopy group\index{stable homotopy group}

\vspace{-3mm}\ 
$$   \pi_{n}^{S}((X_{+}))  = \lim_{ \stackrel{\longrightarrow}{r}}  \  \pi_{n+r}( \Sigma^{r}(X_{+}))  ,   $$
where the limit is taken over the iterated suspension map, depends only on the framed cobordism class of $f$. Incidentally, once one chooses a basepoint in $X$, in the stable homotopy category of  \S\ref{1.3.1} there is a stable homotopy equivalence of the form $X_{+} \simeq  X \vee S^{0}$ and therefore a non-canonical isomorphism $\pi_{*}^{S}(X_{+}) \cong  \pi_{*}^{S}(X) \oplus \pi_{*}^{S}(S^{0})$. 
Pontrjagin (\cite{Pont47,Pont55}) was the first to study stable homotopy groups by means of framed cobordism classes, via the following result, which is proved in \cite[pp.~18--23]{Sto68}:
\begin{theorem}
\label{1.2.6}
The construction of \S\ref{1.2.5} induces an isomorphism of \linebreak graded rings of the form

\vspace{-3mm}\ 
$$   P:     \Omega_{*}^{{\rm fr}}  \stackrel{\cong}{\longrightarrow}  \pi_{*}^{S}(S^{0})  $$
where the ring multiplication in $ \pi_{*}^{S}(S^{0}) $ is given by smash 
product\index{smash product} of maps.
More generally the construction of \S\ref{1.2.5} induces an isomorphism of graded groups of the form

\vspace{-2mm}\ 
$$   P:     \Omega_{*}^{{\rm fr}}(X)  \stackrel{\cong}{\longrightarrow}  \pi_{*}^{S}(X_{+}).   $$
\end{theorem}

In Theorem \ref{1.2.6} the trivial element of $ \pi_{k}^{S}(S^{0}) $ is represented by  the standard $k$-sphere,
$S^{k} \subset {\mathbb R}^{k+1}$ with the framing which is induced by this embedding.

\section{The classical stable homotopy category}\label{categoria}
\label{1.3.1}
In this section we shall give a thumb-nail sketch of the stable homotopy category.  Since we shall only be concerned with classical problems concerning homotopy groups we shall only need the most basic model of a stable homotopy category. 

\smallskip
The notion of a spectrum is originally due to Lima\index{E.L. Lima} \cite{Li58} and was formalised and published in \cite{Wh62}. A spectrum $E$ is a sequence of base-pointed spaces and base point preserving maps (indexed by the positive integers)

\vspace{-5mm}\ 
$$  E :  \   \{  \epsilon_{n} :  \Sigma E_{n}   \longrightarrow  E_{n+1}  \} $$
from the suspension of the $n$-th space $E_{n}$ to the $(n+1)$-th space of $E$. An $\Omega$-spectrum is the same type of data but given in terms of the adjoint maps
$$   E :  \   \{  {\rm adj}(\epsilon_{n}) :   E_{n}   \longrightarrow  \Omega E_{n+1}  \}   $$
from the $n$-th space to the based loops\index{ based loops} on $E_{n+1}$. One requires that the
connectivity\index{connectivity} of the $ {\rm adj}(\epsilon_{n}) $'s increases with $n$. In many examples the ${\rm adj}(\epsilon_{n})$'s are homotopy equivalences which are often the identity map.

\smallskip
According to \cite[p.~140]{Ad74}, a  function $f: E \longrightarrow F$ of degree $r$ in the stable homotopy category is a family of based maps $ \{ f_{n} : E_{n}  \longrightarrow  F_{n-r} \}$ which satisfy the following relations for all $n$:

\vspace{-3mm}\ 
$$    f_{n+1} \cdot \epsilon_{n} =   \phi_{n-r}  \cdot \Sigma f_{n} :  \Sigma E_{n} \longrightarrow  F_{n-r+1} $$
or equivalently
$$   \Omega f_{n+1} \cdot  {\rm adj}(\epsilon_{n}) =  {\rm adj}(\phi_{n-r}) \cdot f_{n} : E_{n} \longrightarrow  \Omega F_{n-r+1}    $$
where $F = \{  \phi_{n} : \Sigma F_{n} \longrightarrow  F_{n+1} \}$.
A morphism $f : E \longrightarrow F$ of degree $r$ is the stable homotopy class of a function $f$.  We shall not need the precise definition of this equivalence relation on functions; suffice to say that that
$[E , F]_{r}$, the stable homotopy classes\index{stable homotopy classes} of morphisms from $E$ to $F$ of degree $r$ in the stable homotopy category, form an abelian group. Summing over all degrees we obtain a graded abelian group $[E , F]_{*}$.

\smallskip
If $X$ and $Y$ are two spaces with base points the smash product $X \wedge Y$ is the product with the two ``axes'' collapsed to a base point. Forming the smash product spectrum $E \wedge F$ of two spectra $E$ and $F$ is technically far from straightforward although a serviceable attempt is made in 
\cite[pp.~158--190]{Ad74}. Essentially the smash product  $E \wedge F$ has an $s$-th space which is constructed from the smash products of the spaces $E_{n} \wedge F_{s-n}$. For our purposes we shall need only the most basic notions of ring spectra and module spectra over ring spectra which are all the obvious translations of the notions of commutative rings and modules over them.  However, one must bear in mind that a commutative ring spectrum\index{commutative ring spectrum} is a generalisation of a graded commutative ring so that the map which involves switching $E_{n} \wedge E_{m}$ to $E_{m} \wedge E_{n}$ carries with it the sign $(-1)^{mn}$ (see \cite[pp.~158--190]{Ad74})!

\begin{example}{\em
\label{1.3.2}\ \vspace{-1ex}

\begin{enumerate}[(i)]
\item The category of topological spaces with base points maps to the stable homotopy category by means of the suspension spectrum\index{suspension spectrum} $\Sigma^{\infty}X$. If $X$ is a based space we define
\[  (\Sigma^{\infty}X)_{0} = X, \  (\Sigma^{\infty}X)_{1} =  \Sigma X ,  \  (\Sigma^{\infty}X)_{2} = \Sigma^{2} X ,  \ldots   \]
with each $\epsilon_{n}$ being the identity map.

\item The fundamental suspension spectrum is the sphere spectrum \linebreak 
$\Sigma^{\infty} S^{0}$ where $S^{0} = \{ 0, 1 \}$ with $0$ as base point. Since $S^{0} \wedge S^{0} = S^{0} $ the suspension spectrum $\Sigma^{\infty} S^{0}$  is a commutative ring spectrum. Since $S^{0} \wedge X = X$ for any based space every spectrum is 
canonically a module spectrum over the sphere spectrum.

\item   If $\Pi$ is an abelian group then the 
Eilenberg-Mac\hspace{1pt}Lane\index{Sammy Eilenberg} 
space\index{Saunders MacLane} \linebreak
$K( \Pi , n)$ is characterised up to homotopy equivalence by the property that
\[  \pi_{i}(   K( \Pi , n) ) \cong  \left\{  \begin{array}{ll}
\Pi   & {\rm if}   \  i=n,  \\
0 & {\rm otherwise} .  
\end{array}  \right.  \]
Therefore  $\Omega K(\Pi , n+1) \simeq K( \Pi , n)$ and we may define the singular homology spectrum\index{singular homology spectrum} or Eilenberg-Mac\hspace{1pt}Lane spectrum\index{Eilenberg-Maclane spectrum} ${\rm H} \Pi$ by  
\[  ({\rm H} \Pi)_{n}  =  K( \Pi , n)  \]
for $n \geq 0$.  We shall be particularly interested in the case of singular homology modulo $2$, ${\rm H}{\mathbb Z}/2$, which is also a commutative ring spectrum, because it is easy to construct non-trivial maps of the form $K( {\mathbb Z}/2 , m) \wedge  K( {\mathbb Z}/2 , n)  \longrightarrow  K( {\mathbb Z}/2 , m+n)$ by means of obstruction theory\index{obstruction theory}.

\item \label{ej4} Every reasonable topological group $G$ has a classifying space $BG$ such that homotopy classes of maps from $W$ to $BG$ correspond to the equivalence classes of principal $G$-bundles over $W$. In particular, this is true for $G= U$, the infinite  unitary group. Homotopy classes of maps from a space $X$ into ${\mathbb Z} \times BU$, whose base point lies in the component $\{ 0 \} \times BU$, classify complex vector bundles\index{complex vector bundle} on $X$ (see \cite{At68} or \cite{Hu66}). The tensor product of the reduced Hopf line bundle\index{Hopf line bundle} on $S^{2}$ with the universal bundle gives a
map $\epsilon : S^{2} \wedge ({\mathbb Z} \times BU) \longrightarrow    {\mathbb Z} \times BU$
and the Bott Periodicity Theorem\index{Bott periodicity Theorem} (\cite{At68} ,  \cite{Hu66}) states that 
the adjoint gives a homotopy equivalence
\[   {\rm adj}(\epsilon) :    {\mathbb Z} \times BU  
\stackrel{\simeq}{\longrightarrow}  \Omega^{2}  ({\mathbb Z} \times BU). \]
This gives rise to the ${\rm KU}$-spectrum or complex periodic K-theory\index{complex periodic K-theory} spectrum defined for $n \geq 0$ by
\[    {\rm KU}_{2n} =   {\mathbb Z} \times BU  \   {\rm and}  \   {\rm KU}_{2n+1} =  \Sigma ( {\mathbb Z} \times BU )  . \]  

Homotopy classes of maps from a space $X$ into ${\mathbb Z} \times BO$ classify real vector bundles\index{real vector bundle} on $X$ (see \cite{At68} or \cite{Hu66}).  The real Bott\index{Raoul Bott} Periodicity Theorem is a homotopy equivalence given by the adjoint of a map, similar to the one in the unitary case,  of the form
\[  \epsilon : S^{8} \wedge ({\mathbb Z} \times BO) \longrightarrow    {\mathbb Z} \times BO  \]
which, in a similar manner, yields a spectrum ${\rm KO}$ whose $8n$-th spaces are each equal to
$  {\mathbb Z} \times BO$.

\item  A connective spectrum\index{connective spectrum} $E$ is one which satisfies
\[  \pi_{r}(E)  =  [  \Sigma^{\infty} S^{r} ,  E   ]_{0}  =  [ \Sigma^{\infty} S^{0} ,  E   ]_{r} = 0   \] 
for $r < n_{0}$ for some integer $n_{0}$. Unitary and orthogonal connective K-theories, denoted respectively by $bu$ and $bo$, are examples of connective spectra with $n_{0}=0$ (\cite[Part III \S16]{Ad74}).

For $m \geq 0$ let $({\mathbb Z} \times BU)(2m, \infty)$ denote a space equipped with a map
\[     ({\mathbb Z} \times BU)(2m, \infty)  \longrightarrow     {\mathbb Z} \times BU  \]
which induces an isomorphism on homotopy groups
\[   \pi_{r}(({\mathbb Z} \times BU)(2m, \infty)) \stackrel{\cong}{ \longrightarrow}  
\pi_{r}(   {\mathbb Z} \times BU) \] 
for all $r \geq 2m$ and such that $\pi_{r}(({\mathbb Z} \times BU)(2m, \infty))=0$ for all $r<2m$. These spaces are constructed using obstruction theory\index{obstruction theory} and are unique up to homotopy equivalence. In particular there is a homotopy equivalence

\vspace{-4mm}\ 
$$    ({\mathbb Z} \times BU)(2m, \infty) \simeq \Omega^{2}   ({\mathbb Z} \times BU)(2m+2, \infty) $$
which yields a spectrum $bu$ given by

\vspace{-4mm}\ 
$$    bu_{2m} =   ({\mathbb Z} \times BU)(2m, \infty)  \  {\rm and }  \   bu_{2m+1} = \Sigma bu_{2m} .$$

The spectrum $bo$ is constructed in a similar manner with 

\vspace{-4mm}\ 
$$  bo_{8m} =   ({\mathbb Z} \times BO)(8m, \infty)   .  $$
In addition one constructs closely related spectra $bso$ and $bspin$ from $BSO = ({\mathbb Z} \times BO)(2, \infty)$ and $BSpin =  ({\mathbb Z} \times BO)(3, \infty)$.

There are canonical maps of spectra $bu \longrightarrow  {\rm KU}$ and 
$bo \longrightarrow  {\rm KO}$ which induce isomorphisms on homotopy groups in dimensions greater than or equal to zero. These are maps of ring spectra.  Similarly we have canonical maps of spectra of the form $bspin \longrightarrow  bso  \longrightarrow  bo$.

\item    The cobordism spectra (\cite[p.~135]{Ad74}) are constructed from Thom spaces (see \cite{Sto68} and Definition \ref{1.2.3}).

For example, let $\xi_{n}$ denote the universal $n$-dimensional vector bundle over $BO(n)$ and let $MO(n)$ be its Thom space\index{Thom space}. The pull-back\index{pull-back} of $\xi_{n+1}$ via the canonical map $BO(n) \longrightarrow  BO(n+1)$ is the vector bundle direct sum $\xi_{n} \oplus 1$
where $1$ denotes the one-dimensional trivial bundle. The Thom space of $\xi_{n} \oplus 1$ is homeomorphic to $\Sigma MO(n)$ which yields a map 
\[  \epsilon_{n} : \Sigma MO(n) \longrightarrow  MO(n+1)  \]
and a resulting spectrum ${\rm MO}$ with  ${\rm MO}_{n} = MO(n)$.

Replacing real vector bundles by complex ones gives ${\rm MU}$ with
\begin{eqnarray*}
& \epsilon_{2n} : \Sigma^{2} MU(n)  \longrightarrow  MU(n+1) = {\rm MU}_{2n+2} &  \\   
& {\rm and}   \quad  {\rm MU}_{2n+1}  =  \Sigma MU(n) .&
\end{eqnarray*}
There are similar constructions associated to the families of classical Lie groups, respectively special orthogonal, spinor and symplectic groups:  ${\rm MSO}$, ${\rm MSpin}$
and ${\rm MSp}$, where ${\rm MSp}_{4n } = MSp(n)$. These cobordism spectra are all connective spectra which are commutative ring spectra by means of the maps (for example, $MO(m) $ $\wedge MO(n)  \longrightarrow  MO(m+n)$) induced by direct sum of matrices in the classical groups.

\item  The following type of non-connective spectrum was introduced in \cite{Sn79}. Suppose that $X$ is a homotopy commutative H-space\index{H-space} and $B \in  \pi_{i}(X)$ or $B \in \pi_{i}^{S}(X) =  \pi_{i}( \Sigma^{\infty} X)$.
Let $X_{+}$ denote the union of $X$ with a disjoint base point. Then there is a stable homotopy equivalence\index{ stable homotopy equivalence} of the form 
$$    \Sigma^{\infty} X_{+}  \simeq   \Sigma^{\infty} X  \vee  \Sigma^{\infty} S^{0} ,  $$

\vspace{-3mm}\noindent 
the wedge sum\index{wedge sum} of the suspension spectra of $X$ and $S^{0}$.
Hence $B \in \pi_{i}(   \Sigma^{\infty} X_{+} )$ and the multiplication in $X$ induces 
$(X \times X)_{+} \cong X_{+} \wedge X_{+}  \longrightarrow  X_{+}$ and thence
$$     \epsilon :  \Sigma^{\infty}  S^{i} \wedge  X_{+}     \stackrel{ B \wedge 1}{\longrightarrow }
\Sigma^{\infty}  X_{+} \wedge X_{+}   \longrightarrow   
\Sigma^{\infty}  X_{+} . $$

\vspace{-1.5mm}\noindent 
This data defines a spectrum $\Sigma^{\infty} X_{+} [1/B]$ in which the $i$-th space is equal to $X_{+}$.
This spectrum will, by construction, be stable homotopy equivalent (via the map $\epsilon$) to its own $i$-th suspension. In particular if $B$ is the generator of $\pi_{2}( {\mathbb CP}^{\infty})$ or $\pi_{2}(BU)$ the resulting spectra $\Sigma^{\infty} {\mathbb CP}^{\infty}_{+} [1/B]$ and $\Sigma^{\infty} BU_{+} [1/B]$ will have this type of stable homotopy periodicity of period $2$.

In Example (\ref{ej4}) we met ${\rm KU}$ which has periodicity of period $2$ and if we add a countable number of copies of ${\rm MU}$ together we may define
\[  {\rm PMU} =   \bigvee_{  n  =  - \infty}^{\infty}  \  \Sigma^{2n} {\rm MU}  \]
which also has periodicity of period $2$. The following result is proved 
in~\cite{Sn79,Sn81}.
\end{enumerate}
}\end{example}

\begin{theorem}
\label{1.3.3}
There are stable homotopy equivalences of the form

\vspace{-3mm}\ 
$$
\Sigma^{\infty} {\mathbb CP}^{\infty}_{+} [1/B]  \simeq  {\rm KU}
\mbox{ \ \ and \ \ }
 \Sigma^{\infty} BU_{+} [1/B]   \simeq  {\rm PMU}  . 
$$
\end{theorem}

\begin{example}{\em
\label{1.3.4}
Here are some ways to make new spectra from old.

\begin{enumerate}[(i)]
\item  \label{1.3.4.1}Given two spectra $E$ and $F$ we can form the smash product spectrum $E \wedge F$. We shall have quite a lot to say about the cases where $E$ and $F$ are various connective K-theory spectra
yielding examples such as $bu \wedge bu$, $bo \wedge bo$ and $bu \wedge bo$. Each of these spectra is a left module spectrum over the left-hand factor and a right module over the right-hand one. For example the left $bu$-module structure on $bu \wedge bo$ is given by the map of spectra

\vspace{-4mm}\ 
$$  bu \wedge ( bu \wedge bo) =  (bu \wedge bu) \wedge bo \stackrel{m  \wedge 1}{\longrightarrow }  
bu \wedge bo $$
where $m : bu \wedge bu \to bu$ is the multiplication in the ring spectrum~$bu$.

\item\label{1.3.4.2} Suppose that $E$ is a spectrum and $G$ is an abelian group. The Moore\index{John Moore} spectrum\index{Moore spectrum} of type 
$G$ (\cite[p.~200]{Ad74},~\cite[Definition~2.1.14]{Rav86}) $MG$ is a connective spectrum characterised by the following conditions on its homotopy and homology groups:
\begin{align*}
&\pi_{r}(MG) =  [ \Sigma^{\infty} S^{0} ,  MG ]_{r}  = 0  \;\,  {\rm \ for } \, 
\  r < 0 ,  \\
&\pi_{0}(MG) \cong  G,  \  {\rm \ and}  \\
&H_{r}(MG ; {\mathbb Z})  =  \pi_{r}(  MG \wedge {\rm H} {\mathbb Z} )   
= 0  \;\,  {\rm \ for }\,  \  r > 0. 
\end{align*}

\vspace{-2mm}The spectrum $E \wedge MG$ is referred to as $E$ with coefficients in $G$\index{$E$ with coefficients in $G$}. For example  $E \wedge \Sigma^{\infty}  \mathbb{RP}^{2}$ is the suspension of 
$E$ with coefficients in ${\mathbb Z}/2$.

\item \label{1.3.4.3} Sometimes one wishes to concentrate on a limited aspect of a spectrum $E$ such as, for example, the $p$-primary part of the homotopy groups. This can often be accomplished using the notions of localisation and completion of spaces or spectra. This appeared first in Alex Zabrodsky's method of ``mixing homotopy types''  \cite{Zab70} which he used to construct non-standard H-spaces.  I first encountered this technique in the mimeographed 1970 MIT notes of Dennis Sullivan\index{Dennis Sullivan} on geometric topology, part of which eventually appeared in \cite{Sul74}. Part I of the original notes is now available as \cite{Sull70}. Localisation of spectra and calculus of fractions is mentioned in\index{Frank Adams} Adams' book \cite{Ad74}
but the final picture was first accomplished correctly by Pete Bousfield\index{Pete Bousfield} in \cite{Bou79}. 

The operations of $p$-localisation and $p$-completion of spectra are related to the spectra with coefficients discussed in (ii). For example, $p$-local and $p$-complete connective K-theories are the subject of \cite{MST77} (see also \cite{Ad78})  and frequently agree with connective K-theories with coefficients 
in ${\mathbb Z}_{(p)}$ (integers localised at $p$) and ${\mathbb Z}_{p}$ (the $p$-adic integers) respectively, when applied to finite spectra.

\item  \label{1.3.4.4}Suppose that $f: E \longrightarrow F$ is a map of spectra of degree zero (if the degree is not zero then suspend or desuspend $F$).  Then there is a mapping cone spectrum\index{mapping cone spectrum} $C_{f}$ called the cofibre of $f$ which sits in a sequence of maps of spectra which generalises the Puppe sequence\index{Puppe sequence} in the homotopy of spaces \cite{Spa66}

\vspace{-8.5mm}\ 
\begin{multline*}
\cdots  \longrightarrow   \Sigma^{-1} C_{f}  \longrightarrow  E  \stackrel{f}{\longrightarrow }  F  \longrightarrow  C_{f}  \\
\longrightarrow  \Sigma E   \stackrel{\Sigma f}{\longrightarrow }  \Sigma F  \longrightarrow  \Sigma C_{f} 
\longrightarrow  \cdots  
\end{multline*}

\vspace{-4mm}\noindent
in which the spectrum to the right of any map is its mapping cone spectrum. This sequence gives rise to long exact homotopy sequences and homology sequences \cite{Ad74} some of which we shall examine in detail later, particularly in relation to $j$-theories and K-theory $e$-invariants\index{e-invariant} defined via Adams operations (as defined in \cite{Ad63}, for example).

 The spectrum $\Sigma^{-1} C_{f} $ will sometimes be called the fibre spectrum\index{fibre spectrum} of $f$.

One may construct new spectra from old as the fibre spectra of maps. For example there is a self-map of
${\rm MU}{\mathbb Z}_{p}$---${\rm MU}$ with $p$-adic coefficients---called the Quillen\index{Dan Quillen} idempotent
\cite{Qu69} (see also \cite[Part II]{Ad74}). The mapping cone of this homotopy idempotent is the Brown-Peterson spectrum
 ${\rm BP}$   (\cite{Rav86},  \cite{Sn02}). There is one such spectrum for each prime. The classical Brown-Peterson spectrum is obtained using coefficients in the $p$-local integers ${\mathbb Z}_{(p)}$ whose $p$-completion is the spectrum obtained using $p$-adic coefficients.

When connective K-theory, say $bu$, is inflicted with coefficients in which $p$ is invertible then there is a self-map $\psi^{p}$ called the $p$-th Adams operation\index{$p$-th Adams operation}. The fibre of 
$\psi^{p}-1$ is an example of a $j$-theory. For example we shall be very interested in 
$$ju  =  {\rm Fibre}(  \psi^{3} - 1 : bu{\mathbb Z}_{2}  \longrightarrow  bu{\mathbb Z}_{2}) . $$

\vspace{-3mm}\noindent
Here ${\mathbb Z}_{p}$ denotes the $p$-adic integers\index{$p$-adic integers}.
We could replace $bu{\mathbb Z}_{2}$ by $2$-localised $bu$ in the sense of  \cite{Bou79}. 

The spectrum $jo$ is defined in a similar manner as (\cite{BJM87}, see also \cite[Chapter Seven]{Arfbk08})
\[   jo  =  {\rm Fibre}(  \psi^{3} - 1 : bo{\mathbb Z}_{2}  \longrightarrow  bspin{\mathbb Z}_{2}) .  \]

\vspace{-2mm}
In \cite{Sn02} one encounters $2$-adic big $J$-theory and $J'$-theory which are defined by the fibre sequences

\vspace{-4mm}\ 
$$
J  \stackrel{\pi}{\longrightarrow}  BP  
\stackrel{\psi^{3} - 1}{\longrightarrow}  BP 
 \stackrel{\pi_{1}}{\longrightarrow} \Sigma J 
$$ 

\vspace{-4mm}
\noindent and
$$J'  \stackrel{\pi'}{\longrightarrow}  BP \wedge BP  
 \stackrel{\psi^{3} \wedge \psi^{3} - 1}{\longrightarrow}  BP \wedge BP 
 \stackrel{\pi'_{1}}{\longrightarrow} \Sigma J'.$$
Here $\psi^{3} $ in $BP$ is induced by the Adams operation $\psi^{3}$ on $MU_{*}(- ;$ $ {\mathbb Z}_{2})$, which commutes with the Quillen idempotent\index{ Quillen idempotent} which is used to define the summand
$BP$.

\item \label{1.3.4.5} One may construct new spectra from old by means of the representability theorem of Ed Brown Jr.~which states that any reasonable generalised cohomology theory (see Definition \ref{1.3.5}) is representable in the stable homotopy category by (i.e. is given by maps into) a unique spectrum. Therefore any procedure, algebraic or otherwise, on a generalised cohomology theory and having sufficient naturality and exactness properties to result in a new generalised cohomology theory results in a new spectrum in the stable homotopy category.

This will be particularly important when we come to the new spectrum which is used in \cite{HHR09}.

The fundamental example is the classical spectrum $BP$ in the category of $p$-local spectra. As we shall see in Theorem \ref{1.3.7} the spectrum $MU$ is intimately related to the universal formal group law. Accordingly when $p$-localised there is an idempotent of $MU$ which corresponds to a result of Cartier concerning the universal $p$-typical formal group law and $BP$ is the resulting stable homotopy analogue. The homotopy of $BP$ is given by $\pi_{*}(BP) = {\mathbb Z}_{(p)}[v_{1}, v_{2},v_{3}, \ldots]$ where $v_{n} \in \pi_{2p^{n}-2}(BP)$ (replace $\mathbb{Z}_{(p)}$ by $\mathbb{Z}_{p}$ to get the homotopy of the $p$-completed $BP$).

The $E_{n}$-spectra were discovered by Jack Morava \cite{Rav86}. There is one for each positive integer $n$ and for each prime $p$ and they are obtained from $BP$. The homotopy $\pi_{*}(E_{n})$ is obtained from $\pi_{*}(BP)$ by
first inverting $v_{n}$ and annihilating the higher-dimensional generators, then completing with respect to the ideal $I_{n} = (p, v_{1}, \ldots ,$ $ v_{n-1})$ and finally adjoining the $(p^{n}-1)$-th roots of unity. For reasons of local class field theory Jonathan Lubin and John Tate \cite{LT66} introduced $\pi_{*}(E_{n})$ and showed that it classifies liftings to Artinian rings of the Honda formal group law $F_{n}$ over ${\mathbb F}_{p^{n}}$.

The $E_{n}$'s are not merely ring spectra. They have a lot of additional structure going by the name of an $E_{\infty}$-ring structure.  A bunch of other $MU$-related spectra may be found in~\cite[p.~124]{Goe04}.

It will come as no surprise to learn that the classical stable homotopy category may be extended to handle spectra with group actions. Brushing the technicalities aside, it should be noted that a result of Mike Hopkins and 
Haynes Miller (\cite[p.~126]{Goe04}) shows that $E_{n}$ is a $G$-spectrum where $G$ is the Morava stabiliser group. In \cite{HHR09} the important motivation for the method comes, I believe, from the study of $G$-spectra acted upon by finite subgroups of Morava stabiliser groups \cite{HHR08}.
\end{enumerate}
}\end{example}

\begin{definition}{\em{\emph{(Generalised homology and cohomology theories)}}
\label{1.3.5}
We shall need the stable homotopy category because it is the correct place in which to study homology and cohomology of spaces and, more generally, spectra.  Generalised homology and cohomology theories originated in \cite{Wh62}.

\smallskip
If $E$ and $X$ are spectra then we define the $E$-homology groups 
of $X$ (\cite[p.~196]{Ad74}) by

\vspace{-5mm}\ 
$$E_{n}(X) =  [ \Sigma^{\infty} S^{0} , E \wedge X ]_{n}$$
and write $E_{*}(X)$ for the graded group given by the direct sum over $n$ of the $  E_{n}(X)$'s.
Sometimes we shall denote $E_{*}(\Sigma^{\infty} S^{0} )$ by $\pi_{*}(E)$, the homotopy of $E$. 

\smallskip
We define the $E$-cohomology groups of $X$ by

\vspace{-3mm}\ 
$$ E^{n}(X) =  [ X , E ]_{-n}  $$
and we write 
$E^{*}(X)$ for the graded group given by the $E$-cohomology groups. 
}\end{definition}

\begin{example}{\em
\label{1.3.6}
Here are some tried and true homology and cohomology theories associated with some of the spectra from Example \ref{1.3.2}.

\begin{enumerate}[(i)]
\item\label{1.3.6.1}  When $E = X = {\rm H}{\mathbb Z}/2$, the mod $2$ Eilenberg-Mac\hspace{1pt}Lane spectrum,
the graded cohomology groups $  ({\rm H}{\mathbb Z}/2)^{*}( {\rm H}{\mathbb Z}/2 )$ form a Hopf algebra 
which is isomorphic to the mod $2$ Steenrod algebra (see \S\ref{1.6.1})
where the product in ${\mathcal A}$ corresponds to composition of maps of spectra. The graded mod $2$ homology groups
$  ({\rm H}{\mathbb Z}/2)_{*}( {\rm H}{\mathbb Z}/2 )$ give the dual Hopf algebra, which is
 isomorphic to the dual Steenrod algebra.
 
\item\label{1.3.6.2} The complex cobordism theories $MU_{*}$ and $MU^{*}$ became very important with the discovery by Daniel Quillen of a connection, described in detail in~\cite[Part II]{Ad74} between formal group laws and $\pi_{*}(MU)$.
\end{enumerate}
}\end{example}

\begin{theorem}{\em(\cite{Mil60}, \cite{Qu69}; see also {\cite[p.~79]{Ad74}})}
\label{1.3.7}\ \vspace{-1ex}

\begin{enumerate}[(i)]
\item \label{1.3.7.1} $\pi_{*}(MU)$ is isomorphic to a graded polynomial algebra\index{graded polynomial algebra} on generators in dimensions $2, 4, 6, 8, \ldots$.

\item\label{1.3.7.2} Let $L$ denote Lazard's ring (\cite[p.~56]{Ad74}) associated to the universal formal group law\index{universal formal group law} over the integers. Then there is a canonical isomorphism
of graded rings

\vspace{-4mm}\ 
$$   \theta : L \stackrel{\cong}{\longrightarrow}  \pi_{*}(MU) .$$
\end{enumerate}
\end{theorem}

Theorem \ref{1.3.7}(i)  was first proved using the classical mod $p$ Adams spectral sequences. Theorem \ref{1.3.7}(ii) is suggested by the fact that the multiplication on $ {\mathbb CP}^{\infty}$ gives rise to a formal group law

\vspace{-6mm}\ 
\begin{multline*}
\pi_{*}(MU)[[x]] =  MU^{*}( {\mathbb CP}^{\infty})  \\ \longrightarrow  MU^{*}( {\mathbb CP}^{\infty} \times  {\mathbb CP}^{\infty})  =  \pi_{*}(MU)[[ x_{1} , x_{2} ]] 
\end{multline*}
with coefficients in $   \pi_{*}( MU)$.

\smallskip
Since Quillen's\index{Dan Quillen} original discovery this theme has developed considerably and details may be found in  \cite{Rav86} (see also  \cite{H07}, \cite{LeMo01}, \cite{LeMo02}, \cite{Lu07}, \cite{Q06}).

\section{Cohomology operations}\label{operaciones}

Every cohomology theory has a ring of stable operations given by its endomorphisms in the stable homotopy category. These algebras of operations are important in the construction of Adams spectral sequences, which are the backbone of the proof in \cite{HHR09} and, indeed, throughout stable homotopy theory \cite{Arfbk08}. Since we should pause to see how Frank Adams constructed the first such spectral sequence  \cite{Ad58} we must begin by recalling the (mod $2$) 
algebra of operations discovered by Norman Steenrod.

\subsection{The Steenrod Operations Modulo $2$}
\label{1.6.1}

Let $H^{*}(X,A ; {\mathbb Z}/2) $ denote singular cohomology\index{singular cohomology} of the pair of topological spaces $A \subseteq X$ (\cite{Spa66}, \cite{StE62}, \cite{Hatch02}). The modulo $2$ Steenrod operations\index{Steenrod operations} are denoted by $Sq^{i}$ for $i \geq 0$ and are characterised by the following axioms (\cite[p.~2]{StE62}):

\begin{enumerate}[(i)]
\item  For all $n \geq 0$, 
$   Sq^{i} :  H^{n}(X,A ; {\mathbb Z}/2)   \longrightarrow  H^{n+i}(X,A ; {\mathbb Z}/2)   $
is a natural homomorphism.

\item $Sq^{0}$ is the identity homomorphism.

\item  If ${\rm deg}(x) = n$ then $Sq^{n}(x) = x^{2}$.

\item If $ i >  {\rm deg}(x)$ then $Sq^{i}(x) = 0$.

\item  The Cartan formula\index{Cartan formula} holds:
$  Sq^{k}(xy) = \sum_{i=0}^{k}  Sq^{i}(x) Sq^{k-i}(y) .  $

\item $Sq^{1}$ is the Bockstein homomorphism associated to the coefficient 
sequence

\vspace{-9mm}\ 
\begin{equation*}
\cdots {\to}\hspace{.2mm} H^{n}(X,A ; {\mathbb Z}/4)
 \hspace{.2mm}{\to}\hspace{.2mm}  H^{n}(X,A ; {\mathbb Z}/2) 
\!\stackrel{Sq^{1}}{ {\to} }  \!H^{n+1}(X,A ; {\mathbb Z}/2) 
\hspace{.2mm}{\to}    \cdots . 
\end{equation*}

\vspace{-3mm}
\item The Adem relations hold\index{Adem relations}. If $0 < a < 2b$ then

\vspace{-4mm}\ 
$$   Sq^{a}Sq^{b} = \sum_{j=0}^{[a/2]}  \   \binom{b-1-j}{a-2j}  \   Sq^{a+b-j} Sq^{j}   $$
where $[y]$ denotes the greatest integer less than or equal to $y$ and 
$
\binom{m}{k} =   \frac{m!}{k! (m-k)!}$  is the usual binomial coefficient (modulo $2$).

\item If $\delta : H^{n}( A ; {\mathbb Z}/2)  \longrightarrow  H^{n+1}(X,A ; {\mathbb Z}/2) $
is the coboundary\index{coboundary} homomorphism from the long exact cohomology sequence of the pair $(X,A)$ then, for all $i \geq 0$,

\vspace{-5mm}\ 
$$  \delta Sq^{i} = Sq^{i} \delta  . $$
\end{enumerate}

\begin{example}{\em{\emph{(Real projective space\index{projective space}  
$\mathbb{RP}^{n}$)}}
\label{1.6.2}
The cohomology ring of real projective $n$-spaces is 
$ H^{*}( \mathbb{RP}^{n} ; {\mathbb Z}/2) \cong  {\mathbb Z}/2[x]/(x^{n+1})$ where ${\rm deg}(x) = 1$. By induction, the axioms of \S\ref{1.6.1} imply 
that (\cite[Lemma 2.4]{StE62})

\vspace{-1mm}\ 
$$  Sq^{i}( x^{k}) =   \binom{k}{i}  \  x^{k+i} .    $$
}\end{example}

\begin{definition}
{\em\emph{(The Steenrod algebra\index{Steenrod algebra} ${\mathcal A}$)}
\label{1.6.3}
\vspace{.5mm}Let $M$ be the graded ${\mathbb F}_{2}$ vector space with basis $\{  Sq^{0}=1, Sq^{1}, Sq^{2}, Sq^{3}, \ldots \}$ with ${\rm deg}(Sq^{i})=i$.  Let $T(M)$ denote the tensor algebra\index{ tensor algebra} of $M$. 
It is a connnected, graded ${\mathbb F}_{2}$-algebra. 

\smallskip
The modulo $2$ Steenrod algebra, denoted by ${\mathcal A}$, is defined to be the quotient of $T(M)$ by the two-sided ideal generated by the Adem relations: if $0 < a < 2b$ then

\vspace{-2mm}\ 
$$    Sq^{a}  \otimes  Sq^{b} = \sum_{j=0}^{[a/2]}  \  \binom{b-1-j}{ a - 2j  } \   Sq^{a+b-j} \otimes  Sq^{j}  .   $$

A finite sequence of non-negative integers $I = (i_{1}, i_{2}, \ldots, i_{k})$ is defined to have length\index{length} $k$, written $l(I) = k$ and moment\index{moment} 

\vspace{-3mm}\ 
$$m(I) = \mmaire\sum_{s=1}^{k} \ si_{s}.$$

\vspace{-.5mm}\noindent 
A sequence $I$ is called admissible\index{admissible} if $\aire i_{1} \geq 1$ and $i_{s-1} \geq 2i_{s}$ for $2 \leq s \leq k$.  Write $Sq^{I} = Sq^{i_{1}} Sq^{i_{2}} \cdots Sq^{i_{k}}$. Then $Sq^{0}$ and all the $Sq^{I}$ with $I$ admissible are called the admissible monomials\index{admissible monomials} of ${\mathcal A}$.
The following result is proved by induction on the moment function $m(I)$.
}\end{definition}

\begin{theorem}{\em (\cite[Theorem 1.3.1]{StE62})}
\label{1.6.4}
The admissible mo\-nomials\index{admissible monomials} form an ${\mathbb F}_{2}$-basis for the Steenrod algebra ${\mathcal A}$.
\end{theorem}

\begin{definition}{\em
\label{1.6.5}
Let $A$ be a graded ${\mathbb F}_{2}$-algebra with a unit $\eta : {\mathbb F}_{2}   \longrightarrow
A$ and a co-unit $\epsilon : A \longrightarrow  {\mathbb F}_{2}$. Therefore $\epsilon \cdot \eta = 1$.
These homomorphisms preserve degree when ${\mathbb F}_{2}$ is placed in degree zero.

\smallskip
Then $A$ is a Hopf algebra\index{Hopf algebra} if: 
\begin{enumerate}[(i)]
\item  There is a comultiplication\index{comultiplication} map
$  \psi :  A \longrightarrow  A \otimes A  $
which is a map of graded algebras when $A \otimes A$ is endowed with the algebra multiplication 
$(a \otimes a') \cdot (b \otimes b') = ab  \otimes a'b'  $ and

\item Identifying  $A \cong A \otimes  {\mathbb F}_{2}   \cong    {\mathbb F}_{2}  \otimes   A$,
$$     1 =  (1 \otimes \epsilon) \cdot  \psi =   ( \epsilon \otimes 1)  \cdot  \psi : A  \longrightarrow  A.   $$
\end{enumerate}
The comultiplication is associative if
$$ ( \psi \otimes 1) \cdot \psi = (1 \otimes \psi) \cdot \psi : A  \longrightarrow  A \otimes A \otimes A.  $$
It is commutative if $\psi =  T \cdot \psi$ where $T(a \otimes a') = a' \otimes a$ for all $a, a' \in A$.
}\end{definition}
\begin{theorem}{\em(\cite[Theorems 2.1.1 and 2.1.2]{StE62})}
\label{1.6.6}
The map of generators

\vspace{-3mm}
$$   \psi( Sq^{k})  = \sum_{i=0}^{k}  Sq^{i} \otimes   Sq^{k-i} $$

\vspace{0mm}\noindent extends to a map of graded ${\mathbb F}_{2}$-algebras
$\psi : {\mathcal A} \longrightarrow  {\mathcal A} \otimes {\mathcal A} $
making ${\mathcal A}$ into a Hopf algebra with a commutative, associative comultiplication.
\end{theorem}

\begin{definition}{\em
\label{1.6.7}
Let ${\mathcal A}^{*}$ denote the dual vector space to ${\mathcal A}$ whose degree $n$ subspace is
${\mathcal A}_{n}^{*} =  {\rm Hom}_{ {\mathbb F}_{2}}( {\mathcal A}_{n},   {\mathbb F}_{2})$.  Let 
\[   \langle - , - \rangle :   {\mathcal A}^{*}  \times {\mathcal A}   \longrightarrow   {\mathbb F}_{2}    \]
denote the evaluation pairing $\langle f , a \rangle = f(a)$.

\smallskip
Let $M_{k} = Sq^{I_{k}}$ where $\mmaire I_{k} = (  2^{k-1}, 2^{k-2}, \ldots , 2, 1) $ for any strictly positive integer $k$. This is an admissible monomial\index{admissible monomial}. Define $\xi_{k} \in {\mathcal A}_{2^{k}-1}^{*}$ by the following formulae for $\langle  \xi_{k} , m  \rangle$ where $m$ runs through all admissible monomials:
\[   \langle \xi_{k} , m  \rangle  =  \left\{  \begin{array}{ll}
1 & {\rm if} \  m = M_{k},  \\
0 & {\rm otherwise}.
\end{array}   \right.   \]
The dual of a commutative coalgebra\index{commutative coalgebra} is a commutative algebra whose multiplication $\psi^{*}$  is the dual of the comultiplication $\psi$. Similarly the dual of a Hopf algebra\index{Hopf algebra} over ${\mathbb F}_{2}$ is again an ${\mathbb F}_{2}$-Hopf algebra.
The following result describes the Hopf algebra ${\mathcal A}^{*}$.
}\end{definition}

\begin{theorem}
{\em(\cite[Theorems 2.2.2 and 2.2.3]{StE62}, \cite{Mil58})}
\label{1.6.8}
The ${\mathbb F}_{2}$-Hopf algebra ${\mathcal A}^{*}$ is isomorphic to the polynomial algebra\index{polynomial algebra} 
${\mathbb F}_{2}[ \xi_{1},$ $ \xi_{2},$ $\xi_{3},$ $ \ldots ,$ $\xi_{k},$  $\ldots ]$
with comultiplication given by 

\vspace{-3mm}\ 
$$   \phi^{*}( \xi_{k})   =   \sum_{i=0}^{k}  \  \xi_{k-i}^{2^{i}} \otimes   \xi_{i} . $$
\end{theorem}
\section{The classical Adams spectral sequence}\label{ASS}
\label{1.4.1}
Now, as promised, here is a sketch of the mod $2$ classical Adams spectral sequence which was first constructed in \cite{Ad58} in order to calculate the stable homotopy groups of spaces. Its generalisation to other cohomology theories can be very technical to construct but this one gives Adams' basic clever idea---namely, to realise by means of maps of spectra the chain complexes which occur in homological algebra. The complete background to the construction in that level of generality is described in the book \cite{MT68}. The construction becomes easier when set up in the stable homotopy category and this is described 
in~\cite[Chapter 3]{Rav86}. The construction of an Adams spectral sequence based on a generalised homology theory was initiated in the case of ${\rm MU}$ by Novikov\index{Sergei Novikov} \cite{Nov67}. These generalised spectral sequences are discussed and explained in \cite{Ad74} and~\cite{Rav86}. 

\subsection{Mod $2$ Adams resolutions}
\label{1.4.2}
Let $X$ be a connective spectrum\index{connective spectrum} such that  
$({\rm H}{\mathbb Z}/2)^{*}( X )\,$---which we shall often write as $H^{*}(X ; {\mathbb Z}/2)\,$---has finite type. This means that each $H^{n}(X ; {\mathbb Z}/2)$ is a finite dimensional ${\mathbb F}_{2}$-vector space. Also $H^{*}(X ; {\mathbb Z}/2)$ is a left module over the mod $2$ Steenrod algebra\index{Steenrod algebra} $\mathcal A$. Therefore one may apply the standard constructions of graded homological algebra\index{graded homological algebra}  (\cite{HS71}, \cite{Mac63})  to form the Ext-groups

\vspace{-3mm}\ 
$$    {\rm Ext}_{{\mathcal A}}^{s,t}(  H^{*}(X ; {\mathbb Z}/2) ,  {\mathbb Z}/2)  .  $$

In order to construct his spectral sequence Frank Adams\index{Frank Adams} imitated the homological algebra construction using spaces---in particular, mod $2$ Eilenberg-Mac\hspace{1pt}Lane\index{Sammy Eilenberg} spaces\index{Saunders MacLane}.

\smallskip
A (mod $2$ classical) Adams resolution\index{Adams resolution} is a diagram of maps of spectra of the following form:
\begin{center}
\begin{picture}(279,62)
\thinlines   
              \put(243,50){$ X_{0}=X$}
              \put(223,52){$ \vector(1,0){10}$}
              \put(225,58){$ g_{0}$}              
              \put(193,50){$ X_{1}$}
              \put(173,52){$ \vector(1,0){10}$}  
              \put(175,58){$ g_{1}$}                           
              \put(143,50){$ X_{2}$}
              \put(123,52){$ \vector(1,0){10}$}
              \put(125,58){$ g_{2}$}               
               \put(93,50){$ X_{3}$}
              \put(73,52){$ \vector(1,0){10}$} 
              \put(75,58){$ g_{3}$}                             
               \put(43,50){$ X_{4}$}
              \put(23,52){$ \vector(1,0){10}$}
              \put(25,58){$ g_{4}$}               
               \put(5,50){$ \cdots$}
              \put(247,43){\vector(0,-1){30}}
              \put(250,28){$f_{0}$}
              \put(197,43){\vector(0,-1){30}}
              \put(200,28){$f_{1}$}
              \put(147,43){\vector(0,-1){30}}
              \put(150,28){$f_{2}$}
              \put(97,43){\vector(0,-1){30}}
              \put(100,28){$f_{3}$}
              \put(47,43){\vector(0,-1){30}}
              \put(50,28){$f_{4}$}
               \put(243,0){$ K_{0}$}
               \put(193,0){$ K_{1}$}
                \put(143,0){$ K_{2}$}
                \put(93,0){$ K_{3}$}
                 \put(43,0){$ K_{4}$}
\end{picture}
\end{center}
in which each $K_{s}$ is a wedge of copies of suspensions of the Eilenberg-Mac\hspace{1pt}Lane spectrum\index{Eilenberg-Mac\hspace{1pt}Lane spectrum} ${\rm H}{\mathbb Z}/2$, each homomorphism    
\[  f_{s}^{*} :    H^{*}(  K_{s} ; {\mathbb Z}/2)  \longrightarrow           H^{*}(  X_{s} ; {\mathbb Z}/2)  \]
is surjective and each 

\vspace{-1mm}\ 
$$X_{s+1}  \stackrel{g_{s}}{\longrightarrow}    X_{s}  \stackrel{f_{s}}{\longrightarrow}   K_{s} $$

\vspace{1mm}\noindent
is a fibring\index{fibring} of spectra (that is, $X_{s+1}$ is the fibre of $f_{s}$).

\smallskip
By Example \ref{1.3.6}(\ref{1.3.6.1}) each $ H^{*}(  K_{s} ; {\mathbb Z}/2)$ is a free ${\mathcal A}$-module and the long exact mod $2$ cohomology sequences split into short exact sequences which splice together to give a free ${\mathcal A}$-resolution of $H^{*}(  X ; {\mathbb Z}/2)$:

\vspace{-5.5mm}\ 
\begin{multline*}
  \cdots  \longrightarrow   H^{*}(  \Sigma^{2} K_{2} ; {\mathbb Z}/2) \longrightarrow H^{*}(  \Sigma K_{1} ; {\mathbb Z}/2)  \\
  \longrightarrow  H^{*}(  K_{0} ; {\mathbb Z}/2) \longrightarrow  H^{*}(  X_{0} ; {\mathbb Z}/2) \longrightarrow   0  .
\end{multline*}
The long exact homotopy sequences of the fibrings of spectra take the form

\vspace{-6mm}\ 
$$   \cdots  \longrightarrow  \pi_{*}(X_{s+1})   \longrightarrow  \pi_{*}(X_{s})  \longrightarrow   \pi_{*}(K_{s})  \longrightarrow  \pi_{*-1}( X_{s+1})    \longrightarrow  \cdots $$
which may be woven together to give an exact couple, which is one of the standard inputs to produce a spectral sequence \cite{Mc01}. In this case the spectral sequence is called the mod $2$ Adams spectral sequence\index{mod $2$ Adams spectral sequence} and satisfies the following properties.
\begin{theorem}
\label{1.4.3}

Let $X$ be a connective spectrum\index{connective spectrum} such that  
$X$ is of finite type. Then there is a convergent spectral sequence\index{convergent spectral sequence} of the form
\[  \begin{array}{c}
 E_{2}^{s,t} = {\rm Ext}_{{\mathcal A}}^{s,t}(  H^{*}(X ; {\mathbb Z}/2) ,  {\mathbb Z}/2 ) 
 \Longrightarrow  \pi_{t-s}(X) \otimes {\mathbb Z}_{2}   
 \end{array}   \]
with differentials\index{differential} $d_{r} : E_{r}^{s,t} \longrightarrow  E_{r}^{s+r,t+r-1}$.
\end{theorem}
\begin{remark}{\em
\label{1.4.4}
In Theorem \ref{1.4.3} the condition that $X$ is of finite type is the spectrum analogue (see \cite{Ad74}) of a space being a CW complex\index{ CW complex} with a finite number of cells in each dimension. It implies that $H^{*}(X ; {\mathbb Z}/2)$ has finite type, but not conversely. Typical useful finite type examples are finite smash products of connective K-theory spectra and suspension spectra of finite CW complexes.
}\end{remark}

\section{Homology operations}\label{homologia}
The approach to non-existence results concerning the Arf-Kervaire invariant used 
in \cite{Akh08a} requires us to know a little about some formulae, one of which is due to J{\o}rgen Tornehave and  me \cite{ST82}, for the invariant in terms of mod $2$ cohomology. In order to be able to sketch this, and several other reformulations of the problem, requires the introduction of the Dyer-Lashof algebra of homology operations---a sort of dual animal to the Steenrod algebra of 
Definition~\ref{1.6.3}.

\begin{definition}{\em
\label{1.7.1}
A based space $X_{0}$ is an infinite loopspace\index{infinite loopspace} if there exists a sequence of based spaces
$X_{1}, X_{2}, X_{3}, \ldots $ such that $X_{i} = \Omega X_{i+1}$, the space of loops in $X_{i+1}$, which begin and end at the base point, for each $i \geq 0$.  A map of infinite loopspaces is defined in the obvious manner.

\smallskip
The principal example of an infinite loopspace is the space 

\vspace{-2mm}
$$QX = \lim\limits_{\stackrel{\longrightarrow}{n}} \Omega^{n} \Sigma^{n}X,$$ the colimit over $n$ of the space of based maps of the $n$-sphere into the $n$-fold suspension of $X$, $\Sigma^{n}X\aire$.

\smallskip
Numerous other examples of infinite loopspaces occur throughout topology---for example, the classifying spaces for topological and algebraic K-theory (see \cite{Ad78}, \cite{MMil79}, \cite{MST77}).

\smallskip
The mod $p$ homology of an infinite loop space admits an algebra of homology operations which complement the operations on homology given by the duals of the Steenrod\index{Norman Earl Steenrod} operations. The homology operations form a Hopf algebra which is usually called the Dyer-Lashof\index{Eldon Dyer} algebra\index{Dick Lashof} after the paper \cite{DL62}. However, the operations originated in the work of Araki\index{S. Araki} and  Kudo\index{T. Kudo} \cite{AK56} (see also \cite{Ni68}). 

\smallskip
We shall be concerned only with homology modulo $2$. Most importantly, the map induced on mod $2$ homology by an infinite loopspace map commutes with the Dyer-Lashof operations.
}\end{definition}

\subsection{The Dyer-Lashof algebra modulo $2$}
\label{1.7.2}

Let $H_{*}( X ; {\mathbb Z}/2)$ denote the singular homology modulo $2$ of an infinite loopspace $X$. 
Hence $H_{*}( X ; {\mathbb Z}/2)$ is an ${\mathbb F}_{2}$-Hopf algebra\index{Hopf algebra}.

\begin{enumerate}[(i)]
\item For each $n \geq 0 $ there is a linear map
\[ Q^{n} : H_{*}(X ; {\mathbb Z}/2) \longrightarrow  H_{*+n}(X ; {\mathbb Z}/2)   \]
which is natural for maps of infinite loopspaces.
\item $Q^{0}$ is the identity map.
\item If ${\rm deg}(x) = n$ then $Q^{n}(x) = x^{2}$.
\item  If ${\rm deg}(x) > n$ then $Q^{n}(x) = 0$.
\item  The Kudo transgression theorem holds:
\[   Q^{n} \sigma_{*} =  \sigma_{*} Q^{n}  \]
where $\sigma_{*} : \tilde{H}_{*}( \Omega X ; {\mathbb Z}/2) \longrightarrow  H_{*+1}(X, {\mathbb Z}/2) $
is the homology suspension map\index{homology suspension map}.
\item The multiplicative Cartan formula holds:
\[   Q^{n}(xy) =  \sum_{r=0}^{n}  \  Q^{r}(x) Q^{n-r}(y) .  \]
\item  The comultiplicative Cartan formula holds:
\[    \psi (Q^{n}(x)) =   \sum_{r=0}^{n}  \  Q^{r}(x') \otimes  Q^{n-r}(x'')   \]
where the comultiplication is given by $\psi(x) = \sum \ x' \otimes x''$.
\item If $\chi :  H_{*}( X ; {\mathbb Z}/2) \longrightarrow  H_{*}( X ; {\mathbb Z}/2)$ is the canonical anti-auto\-mor\-phism\index{canonical anti-automorphism } of the Hopf algebra \cite{MM65} then
\[   \chi \cdot Q^{n} = Q^{n} \cdot \chi  .  \]
\item The Adem relations hold:
\[   Q^{r} \cdot Q^{s} = \sum_{i}  \   \binom{ i-1 }{2i- r}  \   Q^{r+s-i} \cdot Q^{i}   \]
if $r > 2s$.
\item  Let $Sq_{*}^{r} : H_{*}(X ; {\mathbb Z}/2) \longrightarrow  H_{* - i}(X ; {\mathbb Z}/2)$  denote the dual of the Steenrod operation $Sq^{r}$ of \S\ref{1.6.1}. Then the Nishida relations\index{Nishida relations} hold:
\[  Sq_{*}^{r} \cdot Q^{s} = \sum_{i}  \binom{s-r }{r- 2i} \   Q^{s-r+i} \cdot  Sq_{*}^{i} .  \]
\end{enumerate}

\begin{definition}{\em
\label{1.7.3}
$\,$As with Steenrod operations in Definition \ref{1.6.3}, the \linebreak Adem relations for iterated Dyer-Lashof operations lead to the notion of an admissible iterated operation. The element
$Q^{i_{1}} Q^{i_{2}} \cdots Q^{i_{s}}( x)$ is called an admissible\index{admissible} iterated Dyer-Lashof operation if the
sequence $I = (i_{1},$ $i_{2},$ $\ldots ,$ $i_{s})$ satisfies
$i_{j} \leq  2 i_{j+1}$ for $1 \leq j  \leq  s-1$ and 
$i_{u} >   i_{u+1} + i_{u+2} + \ldots + i_{s} + {\rm deg}(x) $ for $1 \leq u \leq s$.  

 The weight\index{weight} of $\mmaire Q^{I} $ is defined to be $2^{s}$.
}\end{definition}

Applications of the Dyer-Lashof operations may be found, for example, 
in~\cite{CLM76,Ko73,MMil79}. Here are a couple.

\subsection{The Kahn-Priddy theorem}
\label{1.5.1}

The motivation for the Kahn-Priddy theorem comes from the classical J-homomorphism \cite{Ad63}, which we shall encounter later (see \S\ref{9.1}) in connection with exotic spheres.

Let $$O_{\infty}({\mathbb R}) = \lim\limits_{\stackrel{\longrightarrow}{m}}  \    O_{m}({\mathbb R})$$ denote the infinite orthogonal group given by the direct limit of the orthogonal groups of $\aire n \times n$ matrices $X$ with real entries and satisfying $X^{-1} = X^{{\rm tr}}$, the transpose of $X$ \cite{Ad73}.  There is a very important, classical homomorphism---the stable J-homomorphism\index{ stable J-homomorphism} \cite{Ad63}---of the form
\[  J :  \pi_{r}(  O_{\infty}({\mathbb R})  )    \longrightarrow    \pi_{r}^{S}(S^{0}) = \pi_{r} ( \Sigma^{\infty} S^{0} ) \]
as a consequence of  the fact that an $n \times n$ orthogonal matrix yields a continuous homeomorphism of $S^{n-1}$. George Whitehead\index{George Whitehead} observed that $J$ factorised through an even more stable J-homomorphism
\[   J' :    \pi_{r}^{S}(  O_{\infty}({\mathbb R})  )   =    \pi_{r}(  \Sigma^{\infty}  O_{\infty}({\mathbb R})  ) 
 \longrightarrow    \pi_{r}^{S}(S^{0}) = \pi_{r} ( \Sigma^{\infty} S^{0} )   \]
and conjectured that $J'$ is surjective when $r > 0$.
Independently,  on the basis of calculations, Mark Mahowald\index{Mark Mahowald} conjectured that $\pi_{r}( \Sigma^{\infty} \mathbb{RP}^{\infty})$ maps surjectively onto the $2$-primary part of $\pi_{r} ( \Sigma^{\infty} S^{0} ) $
via the composition of $J'$ with the well-known map 
\[ \pi_{r}( \Sigma^{\infty} \mathbb{RP}^{\infty})  \longrightarrow     \pi_{r}(  \Sigma^{\infty}  O_{\infty}({\mathbb R})  )  \]
resulting from sending a line in projective space to the orthogonal reflection in its orthogonal hyperplane. Mahowald's conjecture was proved by Dan Kahn\index{Dan Kahn} and Stewart Priddy\index{Stewart Priddy} (\cite{KP72}; see also 
\cite{Ad73}, \cite{KP78} and  \cite{KP78b}). 

\smallskip
The proof uses a construction called the {\em transfer} in stable homotopy together with a careful analysis of its effect when transformed into a map between the infinite loopspaces $Q\mathbb{RP}^{\infty}$ and $Q_{0}S^{0}$, the component of maps of degree zero in $QS^{0}$. The point is that there are isomorphisms
\[   \pi_{*}(\Sigma^{\infty}\mathbb{RP}^{\infty}) \cong \pi_{*}(Q\mathbb{RP}^{\infty})  \    {\rm and}  \  
 \pi_{*}(\Sigma^{\infty}S^{0}) \cong \pi_{*}(Q S^{0}) \]
 in positive dimensions. Mahowald's surjectivity conjecture follows by putting $f$ equal to the J-homomorphism map in the following result. A proof in which the transfer map is replaced by the quadratic part of the Snaith splitting is given in~\cite[\S1.5.10]{Arfbk08}.

\begin{theorem}{\em(The Kahn-Priddy theorem \cite{KP72})}
\label{1.5.10}
Let $f : \mathbb{RP}^{\infty}$ $ \to  Q_{0}S^{0}$ be any map which induces an isomorphism on $\pi_{1}$ and let
$\tilde{f} :  Q \mathbb{RP}^{\infty}  \to  Q_{0}S^{0}$
denote the canonical infinite loopspace map which extends $f$. Then the composite

\vspace{-4mm}\ 
$$   Q_{0}S^{0}    \stackrel{ {\rm transfer} }{\longrightarrow}     Q \mathbb{RP}^{\infty} \stackrel{\tilde{f}}{ \longrightarrow}  Q_{0}S^{0}  $$

\vspace{1mm}\noindent
is a $2$-local homotopy equivalence.
\end{theorem}

\subsection{Hurewicz images in homology modulo $2$}
\label{2.1.1}

Let $H_{*}(W ; {\mathbb Z}/2)$ denote the singular homology modulo $2$ of the space $W$ and let

\vspace{-4mm}\ 
$$   \langle  - , -  \rangle :  H_{*}(W ; {\mathbb Z}/2) \otimes H^{*}(W ; {\mathbb Z}/2) \longrightarrow  {\mathbb Z}/2   $$
denote the canonical non-singular pairing (\cite{Ad74,Hatch02,Spa66}). The dual Steenrod operation\index{dual Steenrod operation} $Sq_{*}^{t}$ is the homology operation\index{homology operation}
\[   Sq_{*}^{t} :  H_{n}(W ; {\mathbb Z}/2)  \longrightarrow  H_{n-t}(W ; {\mathbb Z}/2)  \]
characterised by the equation $\langle  Sq_{*}^{t}(a)  , \alpha \rangle 
= \langle a , Sq^{t}(\alpha)\rangle$ for all $a \in 
 H_{*}(W ; {\mathbb Z}/2) $ and $\alpha  \in   H^{*}(W ; {\mathbb Z}/2) $.

\smallskip Denote by $ H_{*}(W ; {\mathbb Z}/2)_{{\mathcal A}}$ the ${\mathcal A^*}$-annihilated submodule\index{${\mathcal A}$-annihilated submodule}

\vspace{-5mm}\ 
$$   H_{*}(W ; {\mathbb Z}/2)_{{\mathcal A}}  =  \{  x \in  H_{*}(W ; {\mathbb Z}/2)  \  |   \  Sq_{*}^{n}(x) = 0  \   {\rm for  \  all  }  \  n > 0     \} .  $$
Similarly we define the Steenrod annihilated submodule\index{ Steenrod annihilated submodule} $M_{{\mathcal A}}$ for quotients and submodules $M$ of homology which have well-defined actions induced  by the $Sq_{*}^{t}$'s. In particular, $H_{*}(QX ; {\mathbb Z}/2)$ is a Hopf algebra with an action by the $Sq_{*}^{t}$'s and we have $ \underline{Q}H_{*}(QX ; {\mathbb Z}/2)_{{\mathcal A}}$ and 
$ \aire \underline{P}H_{*}(QX ; {\mathbb Z}/2)_{{\mathcal A}}$, the 
${\mathcal A}^*$-annihilated indecomposables and primitives respectively. The latter is of interest because the Hurewicz image of a homotopy class in
$\pi_{r}(QX) \cong \pi_{r}(\Sigma^{\infty}X)$ must lie in $  \underline{P}H_{r}(QX ; {\mathbb Z}/2)_{{\mathcal A}}$ and the primitives are related to the decomposables by a short exact sequence \cite{MM65}.

\smallskip
The following result is proved using the properties listed in \S\ref{1.7.2}. 
\begin{theorem}{\em(\cite{ST82})}
\label{2.1.2}
Let 

\vspace{-2mm}\ 
$$Y = QX = \lim\limits_{\stackrel{\longrightarrow}{m}} \, \Omega^{m} \Sigma^{m}X$$ 

\vspace{-2mm}\noindent
with $X$ a connec\-ted CW complex\index{CW complex}. For $n\geq2$, the 
${\mathcal A}^*$-annihilated decomposables $\aire\underline{Q}H_{*}( Y ; {\mathbb Z}/2)_{{\mathcal A}}$ satisfy:
\begin{enumerate}[(i)]
\item  $   \underline{Q}H_{2^{n}-2}( Y  ; {\mathbb Z}/2)_{{\mathcal A}}   \cong  
H_{2^{n}-2}( X  ; {\mathbb Z}/2)_{{\mathcal A}} $.

\item    $
   \underline{Q}H_{2^{n}-1}( Y  ; {\mathbb Z}/2)_{{\mathcal A}}   \cong  
H_{2^{n}-1}( X  ; {\mathbb Z}/2)_{{\mathcal A}}   \oplus    Q^{2^{n-1}}(H_{2^{n-1}-1}( X  ; {\mathbb Z}/2))$,  
where $Q^{i}\aire$ denotes the $i$-th homology operation\index{homology operation} on $H_{*}( Y  ; {\mathbb Z}/2)$ (\cite{AK56,CLM76,DL62}).
\end{enumerate}
\end{theorem}

\begin{example}{\em
\label{2.1.11}\ \vspace{-1ex}

\begin{enumerate}[(i)]
\item \label{2.1.11.1}If $M$ is a connected closed compact $(2^{n}-2)$-dimensional manifold\index{ manifold} then
$\underline{Q}H_{2^{n}-2}( QM ; {\mathbb Z}/2)_{{\mathcal A}}$ is either zero or ${\mathbb Z}/2$ generated by the fundamental class\index{fundamental class} of $M$, denoted by $[M]$.

\item \label{2.1.11.2}Let $ 0 \not= b_{j} \in  H_{j}( \mathbb{RP}^{\infty} ; {\mathbb Z}/2)$ so that
\[   Sq_{*}^{a}(b_{j}) =    \binom{j-a }{a} \  b_{j-a} . \]
Therefore $\underline{Q}H_{2^{n}-2}( Q  \mathbb{RP}^{\infty}  ; {\mathbb Z}/2)_{{\mathcal A}} = 0$ for $n \geq 2$. Also
\[   \underline{Q}H_{2^{n}-1}( Q  \mathbb{RP}^{\infty}  ; {\mathbb Z}/2)_{{\mathcal A}}  \cong  {\mathbb Z}/2 \langle   b_{2^{n-1}-1}  \rangle  \oplus    {\mathbb Z}/2  \langle  Q^{2^{n-2}}( b_{2^{n-2}-1} )  \rangle  \]
because $\binom{2^{m}-1-a }{a}  \equiv 0$ (modulo $2$).

\item \label{2.1.11.3} Let $O(2)$ denote the orthogonal group\index{orthogonal group} of $2 \times 2$ real matrices $W$ satisfying $W W^{{\rm tr}} = I$, where $W^{{\rm tr}} $ is the transpose\index{transpose}, and let $D_{8}$ denote the dihedral\index{dihedral group} subgroup of order eight generated by $\left(\begin{smallmatrix}0&1\\1&0\end{smallmatrix}\right)$
 and the diagonal matrices. The inclusion of the diagonal matrices gives a chain of groups
\[   {\mathbb Z}/2 \times {\mathbb Z}/2 \subset D_{8} \subset O(2) .  \]
The homology modulo $2$ of these groups is well-known (see~\cite{Pr75,Pr78,Sn89}).
Write $b_{i} * b_{j}$ for the image of $b_{i} \otimes b_{j} \in  H_{i+j}(  \mathbb{RP}^{\infty} \times  \mathbb{RP}^{\infty} ; {\mathbb Z}/2)$ in $H_{i+j}(G ; {\mathbb Z}/2)$ for $G = D_{8}$ or $G = O(2)$. From Theorem \ref{2.1.2} one can show that
\[   H_{2^{n}-2}( BG ; {\mathbb Z}/2)_{{\mathcal A}} =  {\mathbb Z}/2 \langle  b_{2^{n-1}-1} * b_{2^{n-1}-1} \rangle  \]
when $G = D_{8}$ or $G = O(2)$. Hence 
\[  \underline{Q}H_{2^{n}-2}( Q BG ; {\mathbb Z}/2)_{{\mathcal A}} \cong   {\mathbb Z}/2 \langle  b_{2^{n-1}-1} * b_{2^{n-1}-1} \rangle   \]
in these cases.
\end{enumerate}
}\end{example}

\begin{example}{\em(\cite{ST82}, {\cite[Chapter Two]{Arfbk08}})
\label{2.1.14}\ \vspace{-1ex}

\begin{enumerate}[(i)]
\item \label{2.1.14.1}Let $X = \mathbb{RP}^{\infty}$ and let $N_{m}(x_{1}, x_{2}, \ldots )$ denote the $m$-th Newton polynomial\index{Newton polynomial} in variables $x_{1}, x_{2}, \ldots $. From Example \ref{2.1.11}(\ref{2.1.11.2}) and the Milnor-Moore short exact sequence one can show that
\[    \underline{P}H_{2^{n}-2}( Q \mathbb{RP}^{\infty} ; {\mathbb Z}/2)_{{\mathcal A}}   =  {\mathbb Z}/2 \langle   (N_{2^{n-1}-1})^{2}  +  ( Q^{2^{n-2}} N_{2^{n-2}-1} )^{2}  \rangle   .   \]

\item \label{2.1.14.2}From Example \ref{2.1.11}(\ref{2.1.11.3}) one can make similar calculations in the case when $X = BD_{8}$ or $X = BO(2)$. 
\end{enumerate}
}\end{example}

\section{The Arf-Kervaire invariant one problem}\label{AKi1}
\begin{definition}{\em{\emph{(The Arf invariant of a quadratic form)}}
\label{1.8.1}
Let $V$ be a finite dimensional vector space over the field ${\mathbb F}_{2} $ of two elements. A quadratic form\index{ quadratic form} is a function $q : V \longrightarrow   {\mathbb F}_{2}$ such that $q(0) = 0$ and
\[  q(x+y) - q(x) - q(y) = (x , y )  \]
is ${\mathbb F}_{2}$-bilinear (and, of course, symmetric). Notice that $(x , x) = 0$ so that $( - , - )$ is a symplectic bilinear form. Hence ${\rm dim}(V) = 2n$ and to say that $q$ is non-singular
means that there is an ${\mathbb F}_{2}$-basis of $V$,
 $\{a_{1} , \ldots , a_{n}, b_{1}, \ldots , b_{n} \}$ say, such that 
$ (a_{i} , b_{j} )  =  0$ if $i \not= j$, $(a_{i}, b_{i}) = 1$ and $(a_{i}, a_{j} ) = 0 = (b_{i}, b_{j})$ for all $i$ and $j$.

\smallskip
In this case the Arf invariant\index{ Arf invariant} of $q$ is defined to be
\[   c(q) = \sum_{i=1}^{n}  q(a_{i} )   q(b_{i})   \in  {\mathbb F}_{2} .  \]
}\end{definition}

\begin{theorem}{\em(Arf \cite{Arf41}, see 
also~\cite[p.~52]{Brow72} and \cite[p.~340]{Sch85})}
\label{1.8.2}
The invariant $c(q)$ is independent of the choice of basis and two quadratic forms on $V$ are equivalent if and only if their Arf invariants coincide.
\end{theorem}
\subsection{The Arf-Kervaire invariant of a framed manifold}
\label{1.8.3}

Michel Kervaire\index{Michel Kervaire} defined in~\cite{Kerv60} an ${\mathbb F}_{2}$-valued invariant for compact, $(2l-2)$-connnected $(4l-2)$-manifolds 
which are almost parallelisable and smooth in the complement of a point. 
In~\cite{Brow69} Bill Browder extended this definition to any framed, closed $(4l-2)$-manifold where, as in \S\ref{1.2.4},  a framing\index{framing} of $M$ is a homeomorphism
$t:  T(\nu)  \stackrel{\cong}{\longrightarrow}  \Sigma^{N}(M_{+})$, $M_{+}$ is the disjoint union of $M$ and a basepoint and $\nu$ is the normal bundle of an embedding of $M$ into ${\mathbb R}^{4l-2+N}$ and 
$ T(\nu)  =  D(\nu)/S(\nu)$ is the Thom space of the normal bundle. 

\smallskip
Browder\index{Bill Browder} also showed that the Arf invariant of a framed manifold was trivial unless $l = 2^{s}$ for some $s$. He did this by relating the Arf invariant of $M$ to its class in the stable homotopy of spheres and the associated Adams spectral sequence\index{Adams spectral sequence}. 

\smallskip
Recall from Theorem \ref{1.2.6} that the famous Pontrjagin-Thom construction forms the map
\[   S^{4l-2+N} \cong  {\mathbb R}^{4l-2+N}/ ( \infty)  \stackrel{ {\rm collapse}}{\longrightarrow}  
 D(\nu)/S(\nu)  \cong  \Sigma^{N}(M_{+})  \stackrel{ {\rm collapse}}{\longrightarrow}   S^{N}  \]
which  yields an isomorphism between framed cobordism classes of  $(4l-2)$-manifolds and the
 $(4l-2)$-th stable homotopy group of spheres, $\pi_{4l-2}^{S}(S^{0}) =  \pi_{4l-2}( \Sigma^{\infty} S^{0}) \aire$.
 
\smallskip
Here is Bill Browder's\index{Bill Browder} definition, which was simplified by Ed Brown Jr.~in~\cite{EHB72}. Given a framed manifold $M^{2k}$ and $a \in H^{k}(M ; {\mathbb Z}/2) \cong [ M_{+} , K({\mathbb Z}/2 , k) ]$ 
we compose with the Pontrjagin-Thom map to obtain an element of 
\[  \pi_{2k+N}(  \Sigma^{N}K({\mathbb Z}/2 , k) ) \cong {\mathbb F}_{2}.   \]
This is a non-singular quadratic form $q_{M,t}$ on  $H^{k}(M ; {\mathbb Z}/2) $, depending on $t$, and the Arf-Kervaire invariant of $(M,t)$ is $c( q_{M,t} ) \in {\mathbb F}_{2}$.

\begin{example}{\em
\label{1.8.4}
A Lie group has trivial tangent bundle so is frameable. There are framings of $S^{1} \times S^{1}$,
$S^{3} \times S^{3}$ and $S^{7} \times S^{7}$ which have Arf invariant one. As we shall see, there is an elegant way to prove the existence of a framed $M^{30}$ of Arf invariant one (see Corollary \ref{1.8.8}). 
Also in terms of $\pi_{62}^{S}(S^{0})$
a framed manifold of Arf invariant one has been confirmed by the computer calculations of Stan Kochman
\cite{Ko90} and there are also long calculations of Michael Barratt, John Jones and Mark Mahowald \cite{BJM87} asserting the existence in dimension $62$. That is the current extent of existence results.
}\end{example}

\subsection{Equivalent formulations}
\label{1.8.5}
\begin{enumerate}[(i)]
\item \label{1.8.5.1}As above, we write $QS^{0}$ for the infinite loopspace 

\vspace{-4.5mm}\ 
$$ \lim\limits_{ \stackrel{\longrightarrow}{n} } \  \Omega^{n}S^{n},$$

\vspace{-3.5mm}\noindent
the limit over $n\aire$ of the space of based maps from $S^{n}$ to itself. The components of this space are all homotopy equivalent, because it is an H-space, and there is one for each integer. The integer $d$ is the degree of all the maps in the $d$-th component $Q_{d}S^{0}$. The component with $d=1$ is written $SG$, which has an H-space\index{H-space} structure coming from composition of maps. Also $SG$ is important in geometric topology
because $BSG$ classifies stable spherical fibrations and surgery on manifolds starts with the Spivak 
normal bundle, which is a spherical fibration.

Here is a description of the Arf-Kervaire invariant of a framed manifold represented as $\theta \in \pi_{4k+2}( \Sigma^{\infty} S^{0})$ 
(see~\cite[p.~37]{MMil79}). We may form the adjoint, which is a map
$\mmaire{\rm adj}(\theta) : S^{4k+2} \longrightarrow   Q_{0}S^{0}$. Adding a map of degree one yields $ Q_{0}S^{0}
\simeq  SG$ and we may compose with the maps to $SG/SO$ and thence to $G/Top$. These are two of the important spaces which feature in the transformation of geometric surgery theory \cite{Brow72} into homotopy theoretic terms \cite{MMil79}. However, famous work of Dennis Sullivan\index{Dennis Sullivan}, part of his proof of the Hauptvermutung\index{Hauptvermutung}, shows that, at the prime $2$,

\vspace{-3mm}\ 
$$   G/Top \simeq  \prod_{k}  \ K( {\mathbb Z}_{(2)}  , 4k) \times  K({\mathbb Z}/2 , 4k+2) .   $$

\vspace{-3.5mm}\noindent 
Projecting to 
$  K({\mathbb Z}/2 , 4k+2) $ yields an element of $\pi_{4k+2}( K({\mathbb Z}/2 ,$ $4k+2)  )  \cong {\mathbb Z}/2$ which is the Arf-Kervaire invariant of $\theta$. We shall need this description later (see Theorem \ref{1.8.6a}).

\item  \label{1.8.5.2}The Kahn-Priddy Theorem yields a split surjection of the form

\vspace{-2.7mm}\ 
$$   \pi_{m}( \Sigma^{\infty} \mathbb{RP}^{\infty})  \longrightarrow  \pi_{m}( \Sigma^{\infty} S^{0}) \otimes  {\mathbb Z}_{2} $$
for $m>0\,$---where $ {\mathbb Z}_{2} $ denotes the $2$-adic integers\index{ $2$-adic integers}. 
Suppose that \linebreak
$\Theta :  \Sigma^{\infty} S^{2^{n+1}-2} \rightarrow  \Sigma^{\infty} \mathbb{RP}^{\infty}$ is a map whose mapping 
cone is denoted by ${\rm Cone}(\Theta)$. Furthermore (\cite{Brow72})
the image of 

\vspace{-3mm}\ 
$$[ \Theta ] \in  \pi_{2^{n+1}-2}(  \Sigma^{\infty} \mathbb{RP}^{\infty} )$$ under the Kahn-Priddy map has non-trivial Arf-Kervaire invariant if and only if the Steenrod operation (\cite{Sn02,StE62})

\vspace{-6mm}\ 
$$ Sq^{2^{n}} : H^{2^{n} - 1}( {\rm Cone}(\Theta) ; {\mathbb Z}/2)  \cong {\mathbb Z}/2 \rightarrow
H^{2^{n+1} - 1}( {\rm Cone}(\Theta) ; {\mathbb Z}/2)$$
is non-trivial.

Briefly, the reason for this is as follows. 
It is known that the splitting map in the Kahn-Priddy theorem (the transfer) lowers the Adams spectral sequence filtration \cite{K70}. The criterion used in \cite{Brow69} is that an element in the stable homotopy of spheres
has Arf-Kervaire invariant one if and only if it is represented by $h_{2^{n}-1}^{2}$ on the $s=2$ line of the Adams spectral sequence (see Theorem \ref{1.1.2} and Theorem \ref{9.3}). These elements are in filtration two and therefore 
$[ \Theta ] \in  \pi_{2^{n+1}-2}( \Sigma^{\infty} \mathbb{RP}^{\infty} )$ must be in filtration two, one or zero in the Adams spectral sequence for $ \mathbb{RP}^{\infty} $.  One finds that two is impossible and it is easy to show that the filtration cannot be zero, since this would mean that the Hurewicz\index{Witold Hurewicz} image of $[ \Theta ] $ in $H_{*}(\mathbb{RP}^{\infty} ; {\mathbb Z}/2)$ is non-zero. In order to be in filtration one, $[ \Theta ] $ has to be detected by a primary Steenrod operation on the mod $2$ homology of its mapping cone and since the Steenrod algebra is generated by the $\mmaire Sq^{2^{i}}$'s one of these must detect $[ \Theta ] $. The only possibility is $Sq^{2^{n}}\aire$.

A different, more detailed explanation is given in \cite{Ec81}.

\item \label{1.8.5.3} Here are some results from \cite{ST82} (see 
also~\cite[Chapter Two]{Arfbk08}).

By Theorem \ref{1.2.6} a stable homotopy class in $\pi_{2^{n}-2}( \Sigma^{\infty} BO)$ may be considered as a pair $(M , E)$ where $M$ is a frameable $(2^{n}-2)$-manifold and $E$ is a virtual vector bundle on $M$. Theorems \ref{1.8.6} and \ref{1.8.7} use splittings constructed in \cite{Pr78} and \cite{Sn79} together with Dyer-Lashof operations and Steenrod operations in mod $2$ homology of $QX$. 

For $i\geq 1$ let $b_{i} \in H_{i}( \mathbb{RP}^{\infty} ; {\mathbb Z}/2)$ be the generator.
There is an algebra isomorphism
$H_{*}(BO ; {\mathbb Z}/2) \cong {\mathbb Z}/2[ b_{1}, b_{2}, b_{3} , \ldots ]$ where the algebra multiplication is denoted by $x * y$.
\end{enumerate}

The next result follows from Example \ref{2.1.14}(\ref{2.1.14.1}). This is because the adjoint $S^{2^{n+1}-2} \longrightarrow Q_{0}S^{0} \simeq SG$ of a stable homotopy class having non-trivial Arf-Kervaire invariant must have a non-trivial mod $2$ Hurewicz image by \S\ref{1.8.5}(\ref{1.8.5.1}), and the Example \ref{2.1.14}(\ref{2.1.14.1}) shows that there is only one possible non-zero Hurewicz image.
\begin{theorem}
\label{1.8.6a}

Let $\theta \in \pi_{2^{n+1}-2}(\Sigma^{\infty} \mathbb{RP}^{\infty})$ have adjoint 

\vspace{-2.18mm}\ 
$${\rm adj}(\theta) : S^{2^{n+1} -2} \longrightarrow   Q\mathbb{RP}^{\infty}.$$
Let $\lambda \in {\mathbb Z}/2$ denote the Arf-Kervaire invariant of the image of $\theta$ under the Kahn-Priddy map of Theorem \ref{1.5.10}. Then the mod $2$ Hurewicz image of ${\rm adj}(\theta)$ is equal to 

\vspace{-6mm}\ 
$$\lambda \cdot   (N_{2^{n-1}-1})^{2}  +  ( Q^{2^{n-2}} N_{2^{n-2}-1} )^{2}    \in H_{2^{n+1}-2}( Q\mathbb{RP}^{\infty} ; {\mathbb Z}/2) $$
in the notation of Example  \ref{2.1.14}(\ref{2.1.14.1}).
\end{theorem}

The next two results are proved in a similar manner. The formula in Theorem \ref{1.8.7}  will reappear in \S\ref{10.3} where we discuss Pyotr Akhmet'ev's bordism of immersions approach to non-existence of framed manifolds of Arf-Kervaire invariant one (\cite{Akh08a}; see also \cite{Akh08}).

\begin{theorem}
\label{1.8.6}

For $n \geq 3$ there exists a framed $(2^{n}-2)$-manifold with a non-zero Arf-Kervaire invariant
if and only if there exists a stable homotopy class 

\vspace{-6mm}\ 
$$[ M^{2^{n}-2} , E] \in \pi_{2^{n}-2}( \Sigma^{\infty} BO)$$ whose mod $2$ Hurewicz image is equal to 

\vspace{-4mm}\ 
$$b_{2^{n-1}-1} *  b_{2^{n-1}-1} \in H_{2^{n}-2}(BO ; {\mathbb Z}/2).$$ In other words, $b_{2^{n-1}-1} *  b_{2^{n-1}-1}$ is stably spherical.
\end{theorem}

\begin{theorem}
\label{1.8.7}

In Theorem \ref{1.8.6} let ${\rm Arf}[ M^{2^{n}-2} , E] \in {\mathbb Z}/2$ denote the Arf-Kervaire invariant of the framed $(2^{n}-2)$-manifold whose existence is asserted. Then 

\vspace{-4mm}\ 
$$\maire{\rm Arf}[ M^{2^{n}-2} , E] = \langle  [M] ,  w_{2}(E)^{2^{n-1}-1} 
\rangle,$$ where $[M]$ is the mod $2$ fundamental class of $M$ and $w_{2}$ is the second Stiefel-Whitney class.
\end{theorem}

It is very simple to construct, via a balanced product method, a $30$-dimensional manifold $M$ with a $4$-dimensional real vector bundle $E$ over it such that $\langle  [M] ,  w_{2}(E)^{2^{n-1}-1}  \rangle = 1$ in $\mathbb{Z}/2$ (\cite{ST82}; \cite[Chapter~Two]{Arfbk08}). 

\begin{corollary}
\label{1.8.8}
There exists a framed $30$-manifold with non-trivial Arf invariant.
\end{corollary}

The idea of making the requisite framed manifold as a quotient of a balanced product of a surface with $(\mathbb{RP}^{7})^{4} $ was used in~\cite{Jones78} by John Jones, in a highly calculational manner, to prove Corollary~\ref{1.8.8}. 

\subsection{$ju_{*}$-theory reformulation}
\label{1.8.9}

Now let $bu$ denote $2$-adic connective K-theory and define $ju$-theory by means of the fibration

\vspace{-6mm}\ 
$$ju \longrightarrow bu \stackrel{ \psi^{3} - 1}{\longrightarrow}  bu,$$ as in Example \ref{1.3.4}(\ref{1.3.4.4}). Hence $ju_{*}$ is a
generalised homology theory for which  $ju_{2^{n+1}-2}( \mathbb{RP}^{\infty} ) \cong  {\mathbb Z}/2^{n+2}$.
Recall that, if 
%
$$\iota \in ju_{2^{n+1}-2}( S^{2^{n+1}-2}) \cong  {\mathbb Z}_{2}$$
is a choice of generator,
the associated $ju$-theory Hurewicz homomorphism

\vspace{-4mm}\ 
$$ H_{ju} : \pi_{2^{n+1}-2}(  \Sigma^{\infty}   \mathbb{RP}^{\infty} ) \longrightarrow  
ju_{2^{n+1}-2}(  \mathbb{RP}^{\infty} ) \cong  {\mathbb Z}/2^{n+2}  $$

\vspace{1mm}\noindent 
is defined by $H_{ju}([ \theta ] ) = \theta_{*}( \iota )$.

\smallskip
The following result was a conjecture of Barratt-Jones-Mahowald \cite{BJM87}.
\begin{theorem}{\em(\cite{Knapp97}; see also \cite{Sn02}  
and \cite[Chapters Seven and Eight]{Arfbk08})}
\label{1.8.10}
For $n \geq 1$ the image of $[ \Theta ] \in  \pi_{2^{n+1}-2}(  \Sigma^{\infty}  \mathbb{RP}^{\infty} )$ under the
$ju$-theory Hurewicz homomorphism

\vspace{-3mm}\ 
$$H_{ju}([ \Theta ]) \in  ju_{2^{n+1}-2}(  \mathbb{RP}^{\infty} ) \cong  {\mathbb Z}/2^{n+2}$$
is non-trivial if and only if $ Sq^{2^{n}}$ is non-trivial on $ H^{2^{n} - 1}( {\rm Cone}(\Theta) ;  {\mathbb Z}/2) $.

In any case, $2 H_{ju}( [ \Theta] ) = 0$.
\end{theorem}

\section{Exotic spheres, the J-homomorphism and the Arf-Kervaire invariant}
\label{Jmorfismo}
\label{9.1}

The story of the discovery of exotic differentiable structures on spheres is one of the most elegant and important progressive leaps in the history of differential topology---told originally in \cite{KM58}, \cite{KM63}, \cite{Mil56}, \cite{Mil59}; see also \cite{BMil58}. It is worth sketching here because of its relation to the Arf-Kervaire invariant. The theme of the story is that it is not impossibly difficult to make differentiable manifolds which are homotopy equivalent to a sphere, the hard part is to come up with invariants which ensure one has made an exotic sphere.

Let $\Theta_{k}$ denote the group of diffeomorphism classes of smooth manifolds $\Sigma^{k}$ which are homotopy equivalent to $S^{k}$ with group operation induced by connected sum. In 1962 Steve Smale proved the Poincar\'{e} Conjecture in dimensions greater than or equal to five, which implies that $\Sigma^{k}$ is actually homeomorphic to $S^{k}$ in these dimensions. An exotic sphere embeds into Euclidean space with a framing on its normal bundle and, by Theorem \ref{1.2.6}, the Pontrjagin-Thom construction defines an element of $\pi_{k}(\Sigma^{\infty} S^{0})$.

\smallskip
Two framings in the normal bundle of $\Sigma^{k}$ differ by a map into $SO$ so that the above construction yields a homomorphism

\vspace{-3mm}\ 
$$   \tau_{k} : \Theta_{k} \longrightarrow   \pi_{k}(\Sigma^{\infty} S^{0})/{\rm Im}(J)  $$
where $J : \pi_{k}(SO) \longrightarrow  \pi_{k}(\Sigma^{\infty} S^{0})$ is the J-homomorphism introduced in 
\S\ref{1.5.1}. An element in the kernel of $\tau_{k}$ is an exotic sphere which is the boundary of a framed $M^{k+1}$
(see \cite{Mil56} for examples when $k=7$) and a non-trivial element in the image of $\tau_{k}$ comes from an exotic sphere which is not.

\smallskip
Suppose that $M^{n}$ is a framed manifold which is either closed or has an exotic sphere $\Sigma^{n-1}$ as boundary. By surgery (\cite{Brow72}) one may convert $M$, without changing the boundary,  into another framed manifold $W^{n}$, framed cobordant to the $M^{n}$, and which is approximately $n/2$-connected. When $n$ is odd $W^{n}$ will be either a $\Sigma^{n}$ or a disc, whose boundary must be the ordinary sphere $S^{n-1}$. This shows that $\tau_{k}$ is one-one when $k$ is even and onto when $k$ is odd. The obstruction to surjectivity of $\tau_{4l}$ may be shown to vanish using the Hirzebruch signature 
theorem (\cite[\S2]{Mil59}). Furthermore the kernel of $\tau_{4l-1}$ is a large cyclic group which was computed in \cite{KM63}; see also \cite{KM58}. 

\smallskip
The Arf-Kervaire invariant enters the story in order to determine the remaining behaviour of the $\tau_{k}$'s which may be synopsised by the following surgery-based relation

\vspace{-2mm}\ 
$$    {\rm Ker}(\tau_{4l+1}) \oplus {\rm Coker}(\tau_{4l+2}) \cong {\mathbb Z}/2  .  $$

The Arf-Kervaire invariant determines ${\rm Coker}(\tau_{4l+2})$ in the following manner. Suppose that $M^{4l+2}$ is a framed manifold, whose boundary might be an exotic $\Sigma^{4l+1}$. Applying framed surgery we obtain a $2l$-connected framed manifold $\mmaire W^{4l+2}$. The cup-product pairing
\[    H^{2l+1}(W; {\mathbb Z}/2) \otimes  H^{2l+1}(W; {\mathbb Z}/2)  \longrightarrow  H^{4l+2}(W ,  \partial W; {\mathbb Z}/2) \]
evaluates on the fundamental class to give a symmetric, non-singular bilinear form

\vspace{-3mm}\ 
$$   \lambda :  H^{2l+1}(W; {\mathbb Z}/2) \otimes  H^{2l+1}(W; {\mathbb Z}/2)  \longrightarrow {\mathbb Z}/2.   $$
There is a relation $\lambda(x,y) = q_{W,t}(x) + q_{W,t}(y) + q_{W,t}(x+y)  $ where $q_{W,t}$ is as in \S\ref{1.8.3}. This relation, via surgery theory, yields the following result.
\begin{theorem}
\label{9.2}
If there does not exist a framed manifold of dimension $4l+2$ with Arf-Kervaire invariant congruent to $1$ (mod $2$)
then $\tau_{4l+2}$ is surjective and 

\vspace{-5mm}\ 
$$ {\rm Ker}(\tau_{4l+1}) \cong {\mathbb Z}/2.$$
\end{theorem}

In addition we have the following result due to Bill Browder \cite{Brow69} (see also \S\ref{9.4}).
\begin{theorem}
\label{9.3}
A framed $4l+2$ manifold with non-trivial Arf-Kervaire invariant can exist only when $l=2^{j-1}-1$ for some positive integer $j$. If it exists, then it is represented in the classical Adams spectral sequence of Theorem \ref{1.4.3} by
$$h_{j}^{2} \in {\rm Ext}_{{\mathcal A}}^{2, 2^{j+1}}(  {\mathbb Z}/2 , {\mathbb Z}/2 ).$$
\end{theorem}

\subsection{Upper triangular technology (UTT), Browder's Theorem and the Barratt-Jones-Mahowald conjecture}
\label{9.4}

Theorem \ref{9.3} and Theorem \ref{1.8.10} (i.e.~the Barratt-Jones-Mahowald conjecture \cite{BJM87}/Knapp's Theorem \cite{Knapp97}) are related by the UTT technique, which is developed in \cite{Arfbk08}. 

\smallskip
Let $m$ be a positive integer and let $\Theta : \Sigma^{\infty} S^{8m-2} \longrightarrow \Sigma^{\infty}  \mathbb{RP}^{\infty}$ be a morphism in the $2$-local stable homotopy category with mapping cone $C(\Theta)$. From the long exact sequence, there is an isomorphism 

\vspace{-2mm}\ 
$$bu_{8m-1}( C(\Theta)) \cong {\mathbb Z}/2^{4m} \oplus {\mathbb Z}_{2}.$$ Suppose that the Adams operation $\psi^{3}$ on $bu_{8m-1}( C(\Theta))$ satisfies 

\vspace{-2mm}\ 
$$(\psi^{3} - 1)(0,1) = (\frac{3^{4m}-1}{2}, 0),$$ which is equivalent to the $ju$-Hurewicz image of $\Theta$ being non-zero and of order two in $ju_{8m-2}(\mathbb{RP}^{\infty})$. The Barratt-Jones-Mahowald conjecture asserts that this is equivalent to the composition of $\Theta$ with the Kahn-Priddy map having Arf-Kervaire invariant one.

\smallskip
Let $F_{2n}( \Omega^{2} S^{3} )$ denote the $2n$-th filtration of the combinatorial model
for $\Omega^{2} S^{3} \simeq W \times S^{1}$. Let $F_{2n}( W )$
denote the induced filtration on $W$ and let $B(n)$ be the Thom spectrum of the 
canonical bundle induced by $\maire f_{n} :  \Omega^{2} S^{3} \longrightarrow BO$, where $B(0) = S^{0}$ by convention. From \cite{Mah81} one has a left $bu$-module, $2$-local  homotopy equivalence of the form\footnote{In \cite{Arfbk08} and related papers I consistently forgot what I had written in my 1998 McMaster University notes ``On $bu_{*}(BD_{8})$''. Namely, in the description of Mahowald's result I stated that $\Sigma^{4n}B(n) $ was equal to the decomposition factor $F_{4n}/F_{4n-1}$ in the Snaith splitting of $ \Omega^{2} S^{3} $. Although this is rather embarrassing, I got the homology correct so that the results remain correct upon replacing $F_{4n}/F_{4n-1}$ by $\Sigma^{4n}B(n) $ throughout!
I have seen errors like this in
the World Snooker Championship where the no.~1 player misses an easy pot by
concentrating on positioning the cue-ball. In mathematics such errors are
inexcusable whereas in snooker they only cost one the World Championship.}
\[   \bigvee_{n \geq 0 }  bu  \wedge  \Sigma^{4n}B(n)  
\stackrel{ \simeq }{\longrightarrow}  bu  \wedge bo .\]
Therefore $(bu  \wedge bo)_{*}(C(\Theta))  \cong  \bigoplus_{n \geq 0} \  (bu_{*}(C(\Theta) \wedge \Sigma^{4n}B(n) ) $.

\medskip
For $1 \leq k \leq 2m-1$ and $4m \geq 4k - \alpha(k) + 1$ there are isomorphisms of the form (\cite[Chapter Eight \S4]{Arfbk08})

\vspace{-6.3mm}\ 
\begin{align*}
    bu_{8m-1}( C(\Theta)  \wedge \Sigma^{4k}B(k) )   &\cong bu_{8m-1}( \mathbb{RP}^{\infty}  \wedge \Sigma^{4k}B(k) )    \\  &\cong  V_{k} \oplus 
 {\mathbb Z}/2^{4m-4k + \alpha(k) }
\end{align*}

\vspace{-2mm}\noindent 
where $V_{k}$ is a finite-dimensional ${\mathbb F}_{2}$-vector space consisting of elements which are detected in mod $2$ cohomology (i.e. in filtration zero, represented on the $s=0$ line) in the mod $2$ Adams spectral sequence. The map $1 \wedge \psi^{3} \wedge 1$ on $bu \wedge bo \wedge C(\Theta)$ acts on the direct sum decomposition like the upper triangular matrix

\vspace{-2.5mm}\ 
$$\left(\begin{matrix}1 & 1 & 0 & 0 & 0 & \ldots\\0 & 9 & 1 & 0 & 0 & \ldots\\0 & 0 & 9^{2} & 1 & 0 Ê& \ldots\\0 & 0 & 0 & 9^{3} & 1 Ê & \ldots \\\vdots & Ê\vdots & Ê\vdots & Ê\vdots & Ê\vdots & Ê\vdots \end{matrix}\right) . $$

\vspace{0mm}\noindent
In other words $(1 \wedge \psi^{3} \wedge 1)_{*}$ sends the $k$-th summand to itself by multiplication by $9^{k-1}$ and sends the $(k-1)$-th summand to the $(k-2)$-th by a map  

\vspace{-7mm}\ 
$$    \begin{array}{l}
 (\iota_{k,k-1})_{*}   :   V_{k} \oplus 
 {\mathbb Z}/2^{4m-4k + \alpha(k) }   \longrightarrow  V_{k-1} \oplus 
 {\mathbb Z}/2^{4m-4k +4+ \alpha(k-1) }   
 \end{array}   $$

\vspace{-1.5mm}\noindent
for $2 \leq k \leq 2m-1$ and $4m \geq 4k - \alpha(k) + 1$. The right-hand component of this map is injective on the summand $  {\mathbb Z}/2^{4m-4k + \alpha(k) } $ and annihilates $V_{k}$. 
 
From these properties and the formula for $\psi^{3}(0,1)$ on $bu_{8m-1}( C(\Theta))$ one easily obtains a series of equations for the components of $(\eta \wedge 1 \wedge 1)_{*}(0,1)$ where $\eta : S^{0} \longrightarrow bu$ is the unit of $bu$-spectrum. Here we have used the isomorphism $bu_{8m-1}(C(\Theta)) \cong bo_{8m-1}(C(\Theta)) $ since this group is the domain of $(\eta \wedge 1 \wedge 1)_{*}$. This series of equations is impossible unless $m = 2^{q}$ and in that case the $B(2^{q})$-component of $(\eta \wedge 1 \wedge 1)_{*}(0,1)$ is a non-zero class in Adams filtration zero (i.e. a non-zero mod $2$ homology class). This last fact implies that $Sq^{2^{q+2}}$ detects $\Theta$ on its mapping cone. Conversely the system of UTT equations may be used to deduce the formula for $\psi^{3}(0,1)$ from the non-zero homology class in the $B(2^{q})$-component 
(see \cite[Theorems~8.4.6 and~8.4.7]{Arfbk08}).

\smallskip
This upper triangular technology argument simultaneously recovers Theorem \ref{9.3} and  Theorem \ref{1.8.10}.

\smallskip
The UTT technique together with the discussion of \S\ref{9.4} yields the following results. The details are in my unpublished notes. Theorem \ref{9.5}(\ref{9.5.2}) is really a result about lifting Arf-Kervaire invariant one elements to certain types of elements in $ \pi_{2^{n+1}-2}( \Sigma^{\infty} BD_{8})$ because there is a splitting of $\Sigma^{\infty} BD_{8}$ in which 
$\Sigma^{\infty} BPSL_{2}{\mathbb F}_{7}$ is one of the summands~(\cite{MP84}).
\begin{theorem}
\label{9.5}\ \vspace{-1ex}

\begin{enumerate}[(i)]
\item  \label{9.5.1}When $n \geq 5$ there does not exist $\Theta \in \pi_{2^{n+1}-2}( \Sigma^{\infty} \mathbb{RP}^{\infty} \wedge \mathbb{RP}^{\infty} )$ which maps via 

\vspace{-4mm}\ 
$$   \Sigma^{\infty} \mathbb{RP}^{\infty} \wedge \mathbb{RP}^{\infty} \stackrel{H}{ \longrightarrow}   \Sigma^{\infty} \mathbb{RP}^{\infty} \stackrel{KP}{\longrightarrow}   \Sigma^{\infty} S^{0}  $$
to an element having Arf-Kervaire invariant one. Here $H$ is the Hopf construction on the multiplication of $\mathbb{RP}^{\infty} $ and $KP$ is the Kahn-Priddy map.

\item \label{9.5.2}  When $n \geq 5$ there does not exist $\Theta \in \pi_{2^{n+1}-2}( \Sigma^{\infty} BPSL_{2}{\mathbb F}_{7})$  
and a map

\vspace{-7mm}\ 
$$   \Sigma^{\infty} BPSL_{2}{\mathbb F}_{7}    \longrightarrow    \Sigma^{\infty} S^{0}  $$
mapping $\Theta$ to an element having Arf-Kervaire invariant one.
\end{enumerate}
\end{theorem}

\section{Non-existence results for the Arf-Kervaire invariant}
\subsection{Codimension one immersions and the Kervaire invariant one problem}
\label{10.1}\label{noexiste}

The approach to non-existence results for the Arf-Kervaire invariant which is used in \cite{Akh08a} originated in \cite{Ec81} (see \cite{Ec80} for the oriented case).

Let $i : M^{n} \rightarrow {\mathbb R}^{n+1}$ be a self-transverse codimension one immersion of a compact closed smooth manifold. A point of ${\mathbb R}^{n+1}$ is an $r$-fold self-intersection point of the immersion if it is the image under $i$ of at least $r$ distinct points of $M$. Self-transversality implies that the set of $r$-fold intersection points is the image under an immersion of an $(n+1-r)$-dimensional manifold. In particular, the set of  $(n+1)$-fold intersection points is a finite set whose order will be denoted by $\theta(i)$ and whose value mod $2$ is a bordism invariant of the self-transverse immersion.

\smallskip
Suppose that $\xi$ and $\zeta$ are $k$-dimensional real vector bundles. A bundle map $\xi \longrightarrow  \zeta$ is called a $\zeta$-structure on $\xi$. If $\zeta$ is the $k$-plane bundle associated with the universal $G$-bundle where $G \subset O(k)$ then a $\zeta$-structure is traditionally called a $G$-structure. If $\zeta$ is the $k$-bundle over a point this is just a framing. An immersion $i : M^{n} \rightarrow {\mathbb R}^{n+k}$ of a compact closed smooth $n$-manifold will be called a $\zeta$-immersion if its normal bundle has been given a $\zeta$-structure. If such an immersion is self-transverse and $n = km$ then the set of $(m+1)$-fold intersection points is a finite set of order $\theta(i)$ whose value mod $2$ is a bordism invariant of $\zeta$-immersions. The bordism group of $\zeta$-immersions of closed compact $n$-manifolds of codimension $k$ is isomorphic to $\pi_{n+k}( \Sigma^{\infty} T(\zeta))$, where $T(\zeta)$ denotes the Thom space as in Defintion \ref{1.2.3}. This relates bordism of $\zeta$-immersions with the stable homotopy of $\mathbb{RP}^{\infty}$, the Kahn-Priddy theorem and the framed cobordism approach to the stable homotopy groups of spheres.

\smallskip
As for the Arf-Kervaire invariant we have the following result 
from \cite[\S1]{Ec81}.
\begin{theorem}
\label{10.2}

Suppose that $n \equiv 1$ (mod $4$). Then $\theta(i)$ can be odd if and only if there exists a framed $(n+1)$-manifold 
having Arf-Kervaire invariant one.
\end{theorem}
\subsection{Akhmet'ev's approach (\cite{Akh08a}; see also \cite{Akh08})}
\label{10.3}

Suppose that $n = 2^{l}-2$ and that $i: M^{n-1} \rightarrow  {\mathbb R}^{n}$ is a self-transverse codimension one immersion of a compact closed smooth manifold. Let $g : N^{n-2} \rightarrow   {\mathbb R}^{n}$ be the associated 
self-transverse codimension two immersion of the self-intersection submanifold. In Theorem \ref{10.2} the Arf-Kervaire invariant of the framed manifold, whose existence is asserted, is given (\cite[Definition 1]{Akh08a}) by the formula

\vspace{-2mm}\ 
$${\rm Arf}(i) = \langle  [N^{n-2} ] ,  w_{2}(\nu)^{\frac{n-2}{2} } \rangle$$
where $\nu$ is the normal bundle of $g$. This formula is reminiscent of Theorem \ref{1.8.7}, from which it may presumably be deduced.

\smallskip
Armed with this formula the method of  \cite{Akh08a} considers bordism of self-transverse immersions with equivariant $\zeta$-structures. This amounts to studying examples of equivariant stable homotopy groups. The groups involved are the dihedral group of order eight and its products and iterated wreath products. Lifting the problem to these groups is analogous to the Kahn-Priddy theorem which corresponds to lifting to self-transverse immersions with free involution.
Lifting to the dihedral case is presumably analogous to Stewart Priddy's proof that for $i > 0$ $\pi_{i}(\Sigma^{\infty} BD_{8})$ split surjects onto  $\pi_{i}(\Sigma^{\infty} S^{0})$~(\cite{Pr78}). 

\smallskip
Starting from the classification of singular points of generic maps 
$\mathbb{RP}^{s} \rightarrow  {\mathbb R}^{n}$ in the range 
$4s < 3n$~(\cite{Sz97}) Akhmetiev gives a very technical geometric argument to prove the equivariant  self-transverse framed immersions (and skew-framed immersions) cannot exist in high dimensions.

The argument is sufficiently technical to be only accessible to dyed-in-the-wool immersionistos. However, as a start, I recommend the simpler companion paper on the non-existence of framed manifolds having Hopf invariant one \cite{Akh08} as preparatory reading\footnote{Since I finished the first draft of this survey Peter Landweber has
discovered a counterexample to Proposition 41 of~\cite{Akh08}.
In~\cite{La09} he describes the manner in which this error impinges
upon the proofs of the main theorems of~\cite{Akh08} and~\cite{Akh08a}. 
Currently it is unknown to
what extent Proposition 41 is essential to the main results.}\!\!.
\begin{theorem}{\em ({\cite[\S6]{Akh08a}}])}
\label{10.4}
There exists an integer $l_{0}$ such that for all $l \geq l_{0}$ the Arf-Kervaire invariant given by the formula of 
\S \ref{10.3} is trivial.
\end{theorem}

\subsection{The Hill-Hopkins-Ravenel approach \cite{HHR09}}
\label{10.5}

The main theorem of \cite{HHR09} is stronger than that of \cite{Akh08a}  and almost settles the Arf-Kervaire invariant one problem entirely, leaving open only dimension $126$. For a number of reasons this approach has been the more highly publicised. Although the details are very technical the strategy of \cite{HHR09} is quite conceptual because it involves the classical Adams spectral sequence and its generalisations. To get the strategy the odd primary analogue by Doug Ravenel \cite{Rav78} is a good source of preliminary reading. One should bear in mind that there are many Adams-type
spectral sequences, one for each reasonable generalised homology theory. For example, the spectral sequences based on $MU$ and $MSp$ have been used to study $\pi_{*}(\Sigma^{\infty} S^{0})$ but complete determination of these has proved to be too complicated.

\smallskip
If $\tilde{\Omega}$ is a $G$-spectrum for some finite group $G$ there is a spectral sequence of the form

\vspace{-2.9mm}\ 
\[   E_{2}^{s,t} = H^{s}(G ; \pi_{t}( \tilde{\Omega} ) )   \Longrightarrow  \pi_{t-s}( \Omega ) \]
where $\Omega = \tilde{\Omega}^{hG}$, the homotopy fixed point spectrum of $\tilde{\Omega}$. In fact, there exist examples where this spectral sequence {\em is} an Adams spectral sequence
and where $G$ is a Morava stabiliser group of \S\ref{1.3.4}(\ref{1.3.4.4}). This suggests the first step in the strategy---to find a suitable $\tilde{\Omega}$ whose homotopy fixed point spectral sequence receives a map from the Adams-Novikov spectral sequence based on $MU$.

\smallskip
The spectrum $\tilde{\Omega}$ will be an $ E_{\infty}$ spectrum with an action by the cyclic group of order eight $C_{8}$.
Let $MU_{{\mathbb R}}$ denote the spectrum $MU$ with its natural $C_{2}$-involution. Then the four-fold smash product
$MU_{{\mathbb R}} \wedge MU_{{\mathbb R}} \wedge MU_{{\mathbb R}} \wedge MU_{{\mathbb R}}$ has a $C_{8}$-action given by $(a,b,c,d) \mapsto ( \overline{d}, a, b, c)$, which gives a $C_{8}$-spectrum $MU^{((C_{8}))}$. In the equivariant stable homotopy category one may suspend by spheres $S^{V}$, the one-point compactification of a real representation $V$. In the stable category this suspension is invertible which defines suspension by $S^{-V}$. Set

\vspace{-2.5mm}\ 
$$     \tilde{\Omega} = {\rm ho}\hspace{-2.2mm}
\lim\limits_{\stackrel{\longrightarrow}{m}
\;\;\;\;\;} \!\!S^{-ml \rho_{C_{8}}} \wedge MU^{((C_{8}))}  $$

\vspace{-1mm}\noindent
where $l$ is chosen suitably and $\rho_{C_{8}}$ is the real regular representation of $C_{8}$. Then $  \tilde{\Omega} $ is an equivariant $E_{\infty}$ spectrum. 

\smallskip
The relation of $\tilde{\Omega}$ to $MU$ leads to a map of spectral sequences  from the Adams-Novikov spectral sequence

\vspace{-1.5mm}\ 
\[   E_{2}^{s,t} = {\rm Ext}_{MU_{*}(MU)}^{s,t}( MU_{*} , MU_{*}(\Omega))    \Longrightarrow    \pi_{t-s}(\Omega)   \]
to the homotopy fixed point spectral sequence. Furthermore, a key property of $\tilde{\Omega}$ is that its homotopy fixed point spectrum and its actual fixed point spectrum coincide, which makes the latter an easier spectral sequence in which to compute. In particular, it is shown that $\pi_{i}(\Omega)$ vanishes in the range $i = -3, -2, -1$ and that 
$\pi_{i}(\Omega) \cong \pi_{i+256}(\Omega)$ for all $i$ so that $\pi_{2^{n+1}-2}( \Omega) = 0$ for all $n \geq 7$. The choice of the integer $l$ enters in here in order to be able to obtain the correct periodicity.

\smallskip
Modulo lots of tricky technicalities, the argument of \cite{HHR09} concludes as follows. Since, in the Adams-Novikov spectral sequence,  

\vspace{-1mm}\ 
$$   {\rm Ext}_{MU_{*}(MU)}^{0,2^{n+1}-2}  = 0 =  {\rm Ext}_{MU_{*}(MU)}^{1,2^{n+1}-1}$$

\vspace{1mm}\noindent
any element $b_{n}$ in the Adams-Novikov $E_{2}$-term which represents the Arf-Kervaire invariant one element must lie in $\mmaire Ext_{MU_{*}(MU)}^{2,2^{n+1}}$. This is because a map $\aire$ between spectral sequences can only increase the filtration of representatives and, by Theorem \ref{9.3}, the Arf-Kervaire invariant one element is represented by 

\vspace{-2.5mm}\noindent \  
$$h_{n}^{2} \in Ext_{{\mathcal A}}^{2, 2^{n+1}}(  {\mathbb Z}/2 , {\mathbb Z}/2 )$$ 

\vspace{-.5mm}\noindent in the classical Adams spectral sequence. $\aire$The map to the homotopy fixed point spectral sequence faithfully detects all such possible $b_{n}$'s. However $E_{2}^{s,t}= 0$ in the homotopy fixed point spectral sequence for $s \leq 0$ and $t-s$ odd so that for $n \geq 7$ the image of $b_{n}$ cannot be hit by a differential $d_{r}$ with $r \geq 2$. Since $\pi_{2^{n+1}-2}( \Omega) = 0$ the image of $b_{n}$ must be mapped non-trivially by some differential and therefore the same is true for $b_{n}$ in the Adams-Novikov spectral sequence, which implies that $b_{n}$ cannot represent a non-trivial Arf-Kervaire invariant one element. 

\begin{theorem}{\em(\cite{HHR09})}
\label{10.6}
There can exist an element in $\pi_{2^{n+1}-2}(\Sigma^{\infty} S^{0})$
%
%
 having non-trivial Arf-Kervaire invariant only for $n = 1, 2, 3, 4, 5$~and~$6$.
\end{theorem}

This leaves just one unresolved case (see \cite[page (viii)]{Arfbk08}):
\begin{conjecture}
\label{10.7}

There is no element in $\pi_{126}(\Sigma^{\infty} S^{0})$ having Arf-Ker\-vaire invariant equal to one.
\end{conjecture}
\begin{remark}{\em\emph{Having fun!}
\label{10.8}
At this point one of the referees rightly felt entitled to (quote) ``have some fun'' and challenged me to substantiate Conjecture \ref{10.7}. After all the referee's hard work, although I might squirm at the challenge, I have to agree that something should be said. At first sight one might equally well conjecture the opposite, since the methods of \cite{HHR09} break down in dimension $126$. 

\smallskip
In the preface to \cite{Arfbk08} I conjectured that elements $\Theta_{n}$ having non-trivial Arf-Kervaire invariant could only exist for $n = 1, 2, 3, 4, 5$. 

\smallskip
From the 1960's to the 1980's research on the Arf-Kervaire invariant tended to emphasise the constructive aspect. This was very reasonable. I can testify first-hand to the delight one experiences upon constructing a manifold with the desired properties. On the other hand it must have been clear from the 1980's onwards that 
only finitely many $\Theta_{n}$ would exist. This is because there are many equivalent formulations of the Hopf invariant one problem and the Arf-Kervaire invariant one problem which emphasise their similarity and how they form two problems in a systematic progression. In view of the non-existence result for the Hopf invariant one problem (\cite{Ad60})
the only question about Arf-Kervaire invariant one pessimism is where to draw the line.

\smallskip
I imagine that sometime in the 1980's the Arf-Kervaire invariant one practitioners draw the line at what had been constructed by then. Personally, I am aware of a very extensive list of promising Arf-Kervaire invariant one attempts by eminent stable homotopy theorists, which failed above dimension $62$. In addition, myself and others have derived a number of impossibility results such as Theorem \ref{9.5} which, while inconclusive, show that certain constructive methods founder after dimension $62$.

That is enough evidence for me but, of course, it would be great to have a resolution either way of Conjecture \ref{10.7}. 
}\end{remark}

\vspace{-7mm}
\hfill\
{\footnotesize
\parbox{5.15cm}{ $\rule{2mm}{0mm}$ \\ 
$\rule{2mm}{0mm}$ \\ $\rule{2mm}{0mm}$ \\Victor P.~Snaith \\
{\it Department of Pure Mathematics}\\
School of Mathematics and Statistics\\
University of Sheffield\\
Hicks Building, Hounsfield Road\\
Sheffield S3 7RH, UK\\ 
{\sf V.Snaith@sheffield.ac.uk}}} {\hfill}

\label{ultimapagina}
\end{document}